\documentclass[a4paper,11pt,makeidx]{amsart}
\def\mapright#1{\smash{\mathop{\longrightarrow}\limits^{#1}}}
\usepackage{amscd}
\usepackage{xypic}  
\usepackage{amssymb}
\usepackage{amsthm}
\usepackage{epsf}
\makeindex
\newtheoremstyle{break}
  {9pt}
  {9pt}
  {\itshape}
  {}
  {\bfseries}
  {.}
  {\newline}
  {}
\theoremstyle{break}
\newtheorem*{unico}{Theorem}

\theoremstyle{plain}
\newtheorem{thm}{Theorem}[section]
\newtheorem{cor}[thm]{Corollary}

\newtheorem{lemma}[thm]{Lemma}

\newtheorem{prop}[thm]{Proposition}

\newtheorem{defn}[thm]{Definition}
\newtheorem{defn-not}[thm]{Definition-Notation}
\newtheorem{claim}[thm]{Claim}
\newtheorem{rem}[thm]{Remark}
\newtheorem{gerem}[thm]{General remark}

\newtheorem{notation}[thm]{Notation}
\renewcommand{\proofname}{Proof}

\def\Hilb{\operatorname{Hilb}}
\def\Sing{\operatorname{Sing}}
\def\Ker{\operatorname{Ker}}
\def\Coker{\operatorname{Coker}}

\def\min{\operatorname{min}}
\def\Im{\operatorname{Im}}

\def\cork{\operatorname{cork}}
\def\max{\operatorname{max}}
\def\length{\operatorname{length}}
\def\c1{\operatorname{c_1}}
\def\c2{\operatorname{c_2}}
\def\Cliff{\operatorname{Cliff}}
\def\gon{\operatorname{gon}}

\def\CC{{\mathbb C}}
\def\ZZ{{\mathbb Z}}

\def\PP{{\mathbb P}}
\def\Shom{\operatorname{ \mathcal{H}\mathit{om} }}

\def\L{{\mathcal L}}

\def\O{{\mathcal O}}
\def\I{{\mathcal J}}

\def\E{{\mathcal E}}

\def\F{{\mathcal F}}

\def\U{{\mathcal U}}

\def\J{{\mathfrak J}}
\def\*{\otimes}

\def\eqv{\equiv}

\def\+{\oplus}                   
\def\*{\otimes}                  
\def\hpil{\longrightarrow}       
\def\khpil{\rightarrow}

\def\Ext{\operatorname{Ext}}

\def\Pic{\operatorname{Pic}}

\def\Supp{\operatorname{Supp}}

\begin{document}

\title[Surjectivity of Gaussian maps for curves on Enriques surfaces]{Surjectivity of Gaussian maps \\ for 
curves on Enriques surfaces} 

\author[A.L. Knutsen and A.F. Lopez]{Andreas Leopold Knutsen* and Angelo Felice Lopez**}

\address{\hskip -.43cm Andreas Leopold Knutsen, Dipartimento di Matematica, Universit\`a di Roma Tre, Largo San Leonardo
Murialdo 1, 00146, Roma, Italy. e-mail {\tt knutsen@mat.uniroma3.it}}

\thanks{* Research partially supported by a Marie Curie Intra-European Fellowship within the 6th European
Community Framework Programme}

\address{\hskip -.43cm Angelo Felice Lopez, Dipartimento di Matematica, Universit\`a di Roma 
Tre, Largo San Leonardo Murialdo 1, 00146, Roma, Italy. e-mail {\tt lopez@mat.uniroma3.it}}

\thanks{** Research partially supported by the MIUR national project ``Geometria delle variet\`a algebriche" 
COFIN 2002-2004.}


\begin{abstract}
Making suitable generalizations of known results we prove some general facts about Gaussian maps. 
The above are then used, in the second part of the article, to give a set of conditions that insure the
surjectivity of Gaussian maps for curves on Enriques surfaces. To do this we also solve a problem of
independent interest: a tetragonal curve of genus $g \geq 7$ lying on an Enriques surface and general in its linear
system, cannot be, in its canonical embedding, a quadric section of a surface of degree $g - 1$ in $\PP^{g - 1}$.
\end{abstract}

\maketitle

\section{Introduction}
\label{intro}

\noindent Gaussian maps have emerged in the mid 1980's as a useful tool to study the geometry of a given variety $X
\subset \PP^N$ as soon as one has a good knowledge of the hyperplane sections $Y = X \cap H$. 

\noindent Let us briefly recall their definition and notation in the case of curves.

\begin{notation}
{\rm Let $C$ be a smooth irreducible curve and let $L, M$ be two line bundles on $C$. We denote by 
$\mu_{L, M} : H^0(L) \otimes H^0(M) \to H^0(L \otimes M)$ the multiplication map of sections and by $R(L,
M) = \Ker \mu_{L, M}$. The Gaussian map associated to $L$ and $M$ will be denoted by 
\[ \Phi_{L, M} : R(L, M) \to H^0(\omega_C \otimes L \otimes M). \]
This map can be defined locally by $\Phi_{L,M}(s \otimes t) = sdt - tds$ (see \cite{wa}).}
\end{notation}

\noindent Perhaps the first important result, proved by Wahl, who introduced Gaussian maps, is that if a smooth 
curve $C$ lies on a K3 surface, then the Gaussian map $\Phi_{\omega_C, \omega_C}$ cannot be surjective. On the other
hand, as it was proved by Ciliberto, Harris and Miranda \cite{chm}, this map $\Phi_{\omega_C, \omega_C}$ is surjective 
on a curve $C$ with general moduli of genus $10$ or at least $12$. 

\noindent The link with the study of higher dimensional varieties was provided, around the same
period, by Zak, who proved the following result (\cite{za} - see also \cite{bd}, \cite{lv}):

\noindent If $Y \subset \PP^r$ is a smooth variety of codimension at least two with normal
bundle $N_{Y/\PP^r}$ and $h^0(N_{Y/ \PP^r}(-1)) \leq r + 1$, then the only variety $X \subset
\PP^{r + 1}$ that has $Y$ as hyperplane section is a cone over $Y$. 

\noindent Now the point is that, if $Y$ is a curve, we have the formula
\begin{equation*}  
h^0(N_{Y/ \PP^r}(-1)) = r + 1 + \cork \Phi_{H_Y, \omega_Y}
\end{equation*}   
where $H_Y$ is the hyperplane bundle of $Y$. 

\noindent On the other hand, if $Y$ is not a curve one can take successive hyperplane sections of $Y$.
For example, when $X \subset \PP^{r+1}$ is a smooth anticanonically embedded Fano threefold with general
hyperplane section the K3 surface $Y$, in \cite{clm1}, Ciliberto, the second author and Miranda were able to
compute $h^0(N_{Y/ \PP^r}(-1))$ by calculating the coranks of $\Phi_{H_C, \omega_C}$ for the general curve 
section $C$ of $Y$. This then led to recover in \cite{clm1} and \cite{clm2}, in a very simple way, a good 
part of the classification of smooth Fano threefolds \cite{isk1}, \cite{isk2} and of varieties with canonical
curve section \cite{mu}.

\noindent To study other threefolds by means of Zak's theorem, in many cases it is not enough to get down to curve
sections and one needs to bound the cohomology of the normal bundle of surfaces. In \cite{klm} the following general
result was proved:

\begin{prop} 
\label{anysurface}
Let $Y \subset \PP^r$ be a smooth irreducible linearly normal surface and let $H$ be its
hyperplane bundle. Assume there is a base-point free and big line bundle $D_0$ on $Y$
with $H^1(H-D_0)=0$ and such that the general element $D \in |D_0|$ is not rational
and satisfies 
\begin{itemize}
\item[(i)] the Gaussian map $\Phi_{H_D, \omega_{D}(D)}$ is surjective; 
\item[(ii)] the multiplication maps $\mu_{V_D, \omega_D}$ and $\mu_{V_D, \omega_{D}(D)}$ are
surjective, where 

\noindent $V_D := \Im \{H^0 (H-D) \khpil H^0((H-D)_{|D})\}$.
\end{itemize}
Then
\[ h^0(N_{Y/ \PP^r}(-1)) \leq r+1 + \cork \Phi_{H_D, \omega_D}. \]
\end{prop}

\noindent The application of the above proposition clearly points in the following direction: If one
wants to study, with Gaussian maps methods, the existence of threefolds $X \subset \PP^{r+1}$ with given
hyperplane section $Y$, one has to know about the surjectivity of Gaussian maps of type $\Phi_{M, \omega_C}$ 
for curves $C \subset Y$ that are general in their linear system.

\noindent In the present article we do this in the case of Enriques surfaces. 

\noindent This is applied in \cite{klm} to prove the (sectional) genus bound $g \leq 17$ for threefolds
$X \subset \PP^{r+1}$ whose general hyperplane section $Y$ is an Enriques surface. 

\noindent We prove 

\vskip .2cm
 
\begin{unico}
Let $S$ be an Enriques surface and let $L$ be a base-point free line bundle on $S$ with $L^2 \geq 4$. Let $C$ be a general
smooth curve in $|L|$ and let $M$ be a line bundle on $C$. Then the Gaussian map $\Phi_{M, \omega_C}$ is surjective if
one of the hypotheses below is satisfied:
\begin{itemize}
\item[(i)] $L^2 = 4$ and $h^0(4L_{|C} - M) = 0$;
\item[(ii)] $L^2 = 6$ and $h^0((3L + K_S)_{|C} - M) = 0$; 
\item[(iii)] $L^2 \geq 8$ and $h^0(2L_{|C} - M) = 0$;
\item[(iv)] $L^2 \geq 12$ and $h^0(2L_{|C} - M) = 1$;
\item[(v)] $H^1(M) = 0$, $\deg(M) \geq \frac{1}{2} L^2 + 2 \geq 6$ and $h^0(2L_{|C} - M) \leq \Cliff(C) - 2$. 
\end{itemize}
\end{unico}

\vskip .2cm

\noindent The proof of this theorem will be accomplished essentially in two steps. We will first prove, in
Section \ref{basic}, some general facts about Gaussian maps, by generalizing some known results.
Then, in the second step, in Section \ref{onenr}, we will deal with the specific problem of Gaussian maps for
curves on Enriques surfaces. As it turns out, the most difficult point will be to show that a
tetragonal curve of genus $g \geq 7$ lying on an Enriques surface and general in its linear system, in its canonical 
embedding, can never be a quadric section of a surface of degree $g - 1$ in $\PP^{g - 1}$.

\vskip .2cm

\noindent {\it Acknowledgments}. The authors wish to thank Roberto Mu\~{n}oz for several helpful
discussions.

\section{Basic results on Gaussian maps} 
\label{basic}

\noindent We briefly recall the definition, notation and some properties of gonality and Clifford 
index of curves. 

\begin{defn} 
Let $X$ be a smooth surface. We will denote by $\sim$ (respectively $\eqv$) the linear (respectively
numerical) equivalence of divisors (or line bundles) on $X$. We will say that a line bundle
$L$ is {\bf primitive} if $L \eqv kL'$ for some line bundle $L'$ and some integer $k$, implies $k = \pm 1$.
\end{defn} 

\begin{defn} 
Let $C$ be a smooth irreducible curve of genus $g \geq 2$. We denote by $g^r_d$ a
linear system of dimension $r$ and degree $d$ on $C$ and say that $C$ is $k$-{\bf gonal} (and that
$k$ is its {\bf gonality}) if $C$ possesses a $g^1_k$ but no $g^1_{k-1}$. In particular, we call a
$2$-gonal curve {\bf hyperelliptic}, a $3$-gonal curve {\bf trigonal} and a $4$-gonal curve {\bf tetragonal}. 
We denote by $\gon(C)$ the gonality of $C$. 
\end{defn}

\begin{defn} 
Let $C$ be a smooth irreducible curve of genus $g \geq 4$ and let $A$ be a line bundle on
$C$. The {\bf Clifford index} of $A$ is the integer
\[ \Cliff(A) := \deg A - 2(h^0 (A) -1). \]
The {\bf Clifford index of $C$} is
\[ \Cliff(C) := \min \{ \Cliff(A) \hspace{.05in} : \hspace{.05in} h^0 (A) \geq 2, h^1 (A) \geq 2
\}. \]
\noindent We say that a line bundle $A$ on $C$ {\bf contributes to the Clifford index of $C$} if 
$h^0(A) \geq 2, h^1(A) \geq 2$. 
\end{defn}

\subsection{Preliminaries on Gaussian maps}

\noindent We recall some well-known facts about Gaussian maps.

\begin{prop} \cite[Prop.1.10]{wa} 
\label{gaussian} 
Let $C$ be a smooth irreducible nonhyperelliptic curve of genus $g \geq 3$, let $C \subset \PP^{g-1}$
be its canonical embedding and let $M$ be a line bundle on $C$. We have two exact sequences
\begin{equation}
\label{mu}
0 \hpil \Coker \mu_{M, \omega_C} \hpil H^1({\Omega^1_{\PP^{g-1}}}_{|C} \otimes \omega_C \otimes M) \hpil
H^1(M)^{\oplus g} \hpil H^1(\omega_C \otimes M) \hpil 0
\end{equation}
and
\begin{eqnarray}
\label{Phi}
\ \ 0 \hpil \Coker \Phi_{M, \omega_C} \hpil H^1(N_{C / \PP^{g-1}}^{\ast} \otimes \omega_C \otimes
M) & \hpil & H^1({\Omega^1_{\PP^{g-1}}}_{|C} \otimes \omega_C \otimes M) \hpil \\
\nonumber & \hpil & H^1(\omega_C^2 \otimes M) \hpil 0.
\end{eqnarray}
In particular
\begin{itemize}
\item[(a)] if $H^0(N_{C / \PP^{g-1}} \otimes M^{-1} ) = 0$ then $\Phi_{M, \omega_C}$ is surjective;
\item[(b)] if $H^1(M) = 0$ and $\mu_{M, \omega_C}$ is surjective then 
$\cork \Phi_{M, \omega_C} = h^0(N_{C / \PP^{g-1}} \otimes M^{-1})$.
\end{itemize}
\end{prop}

\noindent In the sequel we will collect some results about Gaussian maps of type
$\Phi_{M, \omega_C}$ for curves $C$ of low genus or low gonality or with Clifford index higher than
$h^0(2K_C - M) + 2$.

\noindent We start with an elementary but useful fact.

\begin{lemma}
\label{quadrics}
For $a \geq 2$, let $Q_1, \ldots, Q_a$ be linearly independent homogeneous polynomials
of degree 2 in $X_0, \ldots, X_r$. Suppose that the relations among $Q_1, \ldots, Q_a$ are
generated by $R_i = [R_{i1}, \ldots, R_{ia}]$, for $1 \leq i \leq b$. If $(c_1, \ldots, c_a) \in 
\CC^a - \{0\}$ then there exists an $i$ such that 

\[ \sum\limits_{j=1}^a c_j R_{ij} \neq 0. \]
\end{lemma}
\noindent {\it Proof}.
Suppose to the contrary that $\sum\limits_{j=1}^a c_j R_{ij} = 0$ for every $i$ with $1 \leq i \leq b$.

\noindent Without loss of generality assume that $c_1 \neq 0$, so that 
\begin{equation}
\label{eqr}
R_{i1} = - \sum\limits_{j=2}^a c_1^{-1} c_j R_{ij}, \ \ 1 \leq i \leq b.
\end{equation}

\begin{claim}
\label{qprimo}
Set $Q_1' = Q_1, Q_j' = Q_j - c_1^{-1} c_j Q_1$ for $2 \leq j \leq a$. Then
\begin{itemize}
\item[(i)] $Q_1', \ldots, Q_a'$ are linearly independent;
\item[(ii)] the relations among $Q_1', \ldots, Q_a'$ are generated by $S_i = [0, R_{i2}, \ldots,
R_{ia}]$, for $1 \leq i \leq b$.
\end{itemize}
\end{claim}

\begin{proof}
Consider a relation $\sum\limits_{j=1}^a R_j' Q_j' = 0$, where the $R_j'$'s are polynomials. Then 

\noindent $R_1' Q_1 + \sum\limits_{j=2}^a R_j' (Q_j - c_1^{-1} c_j Q_1) = 0$, whence 
\begin{equation}
\label{eqr1}
(R_1' - \sum\limits_{j=2}^a c_1^{-1} c_j R_j') Q_1 + \sum\limits_{j=2}^a R_j' Q_j = 0. 
\end{equation}

\noindent If all $R_j'$'s are complex numbers we get $R_j' = 0$ for all $j$, proving (i).

\noindent To see (ii), by (\ref{eqr1}) and the hypothesis of the lemma we deduce that there are polynomials 
$d_j$ such that
\[ [R_1' - \sum\limits_{j=2}^a c_1^{-1} c_j R_j', R_2', \ldots, R_a'] = \sum\limits_{i=1}^b d_i R_i =
[\sum\limits_{i=1}^b d_i R_{i1}, \sum\limits_{i=1}^b d_i R_{i2}, \ldots, \sum\limits_{i=1}^b d_i R_{ia}] \]
whence $R_j' = \sum\limits_{i=1}^b d_i R_{ij}$ for $2 \leq j \leq a$ and
\[R_1' = \sum\limits_{j=2}^a c_1^{-1} c_j R_j' + \sum\limits_{i=1}^b d_i R_{i1} = \sum\limits_{i=1}^b d_i
(\sum\limits_{j=2}^a c_1^{-1} c_j R_{ij} + R_{i1}) = 0 \] 
by (\ref{eqr}). Now
\[\sum\limits_{i=1}^b d_i S_i = [0, \sum\limits_{i=1}^b d_i R_{i2}, \ldots , \sum\limits_{i=1}^b d_i R_{ia}]
= [R_1', R_2', \ldots, R_a']. \]
\end{proof}

\renewcommand{\proofname}{Conclusion of the proof of Lemma {\rm \ref{quadrics}}}  
\begin{proof}
Consider the Koszul relation $[Q_2', - Q_1', 0, \ldots, 0]$ among $Q_1', \ldots, Q_a'$. By the claim there are 
polynomials $d_i$ such that $\sum\limits_{i=1}^b d_i S_i = [Q_2', - Q_1', 0, \ldots, 0]$, giving the
contradiction $Q_2' = 0$.
\end{proof}
\renewcommand{\proofname}{Proof}

\noindent In many cases, to compute the corank of Gaussian maps, or, as in Proposition \ref{gaussian}, to compute a
suitable cohomology group involving the normal bundle, it is quite convenient to know some surface
containing the given curve. The result below will help to compute the cohomology of the normal bundle with the
help of the  surface.

\begin{lemma}
\label{N2}
Let $Y \subset \PP^r$ be an integral subvariety that is scheme-theoretically intersection of quadrics and let
$X \subset Y$ be a smooth irreducible nondegenerate subvariety. Let $L = \O_Y(1)$ and $M$ a line bundle on $X$. 
Suppose that either
\begin{itemize}
\item[(i)] $h^0(2L_{|X} - M) = 0$ 
or
\item[(ii)] $h^0(2L_{|X} - M) = 1$ and the relations among the quadrics cutting out $Y$ are generated by 
linear ones.
\end{itemize}
Let $\F_{X, Y} = \Shom_{\O_{\PP^r}}(\I_{Y / \PP^r}, \O_X)$. Then $H^0(\F_{X, Y} \otimes M^{-1}) = 0$.
\end{lemma}

\begin{rem}
\label{N2rem} 
{\rm When $Y$ is smooth we have that $\F_{X, Y} = {N_{Y / \PP^r}}_{|X}$. The fact that $Y \subset
\PP^r$ is scheme-theoretically intersection of quadrics certainly holds if $Y$ satisfies property $N_1$, that 
is $Y$ is projectively normal and its homogeneous ideal is generated by quadrics (\cite[Def.1.2.5]{la1},
\cite{gr}). Also the fact that the relations among the quadrics cutting out $Y$ are generated by linear ones
certainly holds if $Y$ satisfies property $N_2$, that is $Y$ satisfies property $N_1$ and the relations among
the quadrics generating its homogeneous ideal are generated by linear ones (\cite[Def.1.2.5]{la1}, \cite{gr}).
The difference, in our case, is that we do not assume $Y$ to be linearly normal.}
\end{rem} 
 
\renewcommand{\proofname}{Proof of Lemma {\rm \ref{N2}}}
\begin{proof}  
Let $\{Q_1, \ldots, Q_a \}$ be linearly independent quadrics cutting out $Y$ scheme-theoretically and consider the
corresponding beginning of the minimal  free resolution of $\I_{Y / \PP^r}$:
\[ \bigoplus_{i \geq 0} \O_{\PP^r} (-3 - i)^{\oplus b_i} \hpil \O_{\PP^r} (-2)^{\oplus a} \hpil \I_{Y / \PP^r} 
\hpil 0. \]  
Applying the left exact functor $\Shom_{\O_{\PP^r}}(- , \O_X)$ we get an exact sequence
\[ 0 \hpil \F_{X, Y} \hpil \O_X (2)^{\oplus a} \hpil \bigoplus_{i \geq 0} \O_X (3 + i)^{\oplus b_i} \]
whence an exact sequence
\[ 0 \hpil H^0(\F_{X, Y} \otimes M^{-1}) \hpil H^0(2L_{|X} - M)^{\oplus a} \mapright{\varphi} \bigoplus_{i
\geq 0} H^0((3+i)L_{|X} - M)^{\oplus b_i}. \]
Then $H^0(\F_{X, Y} \otimes M^{-1}) = \Ker \varphi$.

\noindent If we are under hypothesis (i), then obviously $H^0(\F_{X, Y} \otimes M^{-1}) = 0$.

\noindent If we are under hypothesis (ii), then $b_i = 0$ for $i \geq 1$ and we will prove that $\Ker \varphi = 0.$ 

\noindent To this end let $\sigma$ be a generator of $H^0(2L_{|X} - M)$. For $1 \leq i \leq b_0$ let 
$R_i = [R_{i1}, \ldots, R_{ia}]$ be the linear relations generating all relations among $Q_1, \ldots, Q_a$, so 
that the map $\varphi$ is given by the matrix $({R_{ij}}_{|X})$. If $0 \neq (c_1 \sigma, \ldots, c_a \sigma) \in
\Ker \varphi$ then, for every $i$ such that $1 \leq i \leq b_0$, we have $\sum\limits_{j=1}^a {R_{ij}}_{|X} c_j
\sigma = 0$ whence $(\sum\limits_{j=1}^a c_j R_{ij} )_{|X} = 0$. As $X$ is nondegenerate and
$\sum\limits_{j=1}^a c_j R_{ij}$ is a linear polynomial, we deduce that $\sum\limits_{j=1}^a c_j R_{ij} = 0$ for all $i$
with $1 \leq i \leq b_0$, contradicting Lemma \ref{quadrics}.
\end{proof}
\renewcommand{\proofname}{Proof}

\noindent Now the first general result about Gaussian maps. 

\begin{prop}
\label{low}
Let $C$ be a smooth irreducible nonhyperelliptic curve of genus $g \geq 3$ and let $M$ be a line bundle on
$C$. We have
\begin{itemize}
\item[(a)] If $g = 3$ then $\cork \Phi_{M, \omega_C} \geq h^0(4K_C - M) - \cork \mu_{M, \omega_C} - 3
h^1(M)$, with equality if $H^0(- M) = 0$.
\item[(b)] If $g = 4$ then $\cork \Phi_{M, \omega_C} \geq h^0(2K_C - M) + h^0(3K_C - M) -
\cork \mu_{M, \omega_C} - 4 h^1(M)$, with equality if $H^0(- M) = 0$.
\item[(c)] If $g = 5$ and $C$ is nontrigonal then $\cork \Phi_{M, \omega_C} \geq 3 h^0(2K_C - M) - \cork
\mu_{M, \omega_C} - 5 h^1(M)$, with equality if $H^0(- M) = 0$.
\item[(d)] Suppose that $C$ is a plane quintic and $A$ is the very ample $g^2_5$ on $C$. If $H^0(5A - M) = 0$ 
then $\Phi_{M, \omega_C}$ is surjective. If $H^1(M) = 0$ and $\mu_{M, \omega_C}$ is surjective then $\cork
\Phi_{M, \omega_C} \geq h^0(5A - M)$, with equality holding if in  addition $h^0(4A - M) \leq 1$.
\item[(e)] Suppose that $C$ is trigonal, $g \geq 5$ and $A$ is a $g^1_3$ on $C$. If $h^0(2K_C - M) \leq
1$ and $H^0(3K_C - (g - 4) A - M) = 0$ then $\Phi_{M, \omega_C}$ is surjective. If $H^1(M) = 0$
and $\mu_{M, \omega_C}$ is surjective then $\cork \Phi_{M, \omega_C} \geq h^0(3K_C - (g - 4) A - M)$, with
equality holding if in addition $h^0(2K_C - M) \leq 1$.
\end{itemize}
\end{prop}

\begin{proof}
Assertions (a), (b) and (c) follow easily from Proposition \ref{gaussian}.

\noindent Let us prove (d). In the canonical embedding $C \subset \PP^5$ we have that $C$ is contained in the
Veronese surface $Y$ and we have an exact sequence
\begin{equation}
\label{ver}
0 \hpil N_{C / Y} \otimes M^{-1} \hpil N_{C / \PP^5} \otimes M^{-1} \hpil {N_{Y / \PP^5}}_{|C} \otimes M^{-1}
\hpil 0.
\end{equation}
Observe that $h^0(N_{C / Y} \otimes M^{-1}) = h^0(5A - M)$. Now if $h^0(5A - M) = 0$ then also $h^0(2K_C - M) = 
h^0(4A - M) = 0$ and from (\ref{ver}) and Proposition \ref{gaussian} (a), we see that to prove (d) we just need to
show that $H^0({N_{Y / \PP^5}}_{|C} \otimes M^{-1}) = 0$. The latter follows by Lemma \ref{N2} and Remark \ref{N2rem}
since, as is well-known, $Y$ satisfies property $N_3$.

\noindent Now if $H^1(M) = 0$ and $\mu_{M, \omega_C}$ is surjective, we have that $\cork \Phi_{M, \omega_C} = 
h^0(N_{C / \PP^5} \otimes M^{-1}) \geq h^0(5A - M)$ by Proposition \ref{gaussian} (b) and (\ref{ver}). If we also assume
that $h^0(4A - M) = h^0(2K_C - M) \leq 1$ then we can apply again Lemma \ref{N2} and Remark \ref{N2rem}. We get that
$h^0({N_{Y / \PP^5}}_{|C} \otimes M^{-1}) = 0$, whence, from (\ref{ver}), that $h^0(N_{C / \PP^5} \otimes M^{-1}) =
h^0(5A - M)$.

\noindent To see (e) recall that, in the canonical embedding $C \subset \PP^{g - 1}$, we have
\cite[6.1]{sch1} that $C \in |3H - (g - 4) R|$ on a rational normal surface $Y \subset \PP^{g - 1}$, where $H$
is its hyperplane bundle and $R$ its ruling. Since, as is well-known, $Y$ satisfies property $N_{g - 3}$,
applying, as in case (d), Lemma \ref{N2} and Proposition \ref{gaussian} we get (e).
\end{proof}

\noindent Note that the cases (d), (e) of the above proposition and the corollary below are a slight
improvement of \cite[Thm.2.4]{te} (because we also consider the case $h^0(2K_C - M) = 1$).

\begin{cor} 
Let $C$ be a smooth irreducible curve of genus $g \geq 5$ and let $M$ be a line bundle on $C$.
Then the Gaussian map $\Phi_{M, \omega_C}$ is surjective if one of the hypotheses below is satisfied:
\begin{itemize}
\item[(a)] $C$ is a plane quintic and $\deg M \geq 25$, $M \neq 5A$ if equality holds, where $A$ is the very
ample $g^2_5$ on $C$;
\item[(b)] $C$ is trigonal and $\deg M \geq \max \{4g - 6, 3g + 6 \}$, $M \neq 3K_C - (g - 4) A$ if
$g \leq 12$ and $\deg M = 3g + 6$, where $A$ is a $g^1_3$ on $C$.
\end{itemize}
\end{cor}

\begin{proof}
(a) follows immediately from Proposition \ref{low}(d) while (b) is a consequence of Proposition
\ref{low}(e) since, if $h^0(2K_C - M) \geq 2$, then $\deg (2K_C - M) \geq 3$, a contradiction.
\end{proof}

\noindent Another easy but useful consequence of the proof of Lemma \ref{N2} is the following.

\begin{prop}
\label{Cliff}
Let $C$ be a smooth irreducible curve of genus $g \geq 5$ and let $M$ be a line bundle on $C$. Suppose 
that either
\begin{itemize}
\item[(i)] $\Cliff(C) = 2$ and $h^0(2K_C - M) = 0$ 
or
\item[(ii)] $\Cliff(C) \geq 3$ and $h^0(2K_C - M) \leq 1$.
\end{itemize}
Then $\Phi_{M, \omega_C}$ is surjective.
\end{prop}

\begin{proof}
Since $\Cliff(C) \geq 2$, by \cite{vo}, \cite{sch2}, the resolution of the ideal sheaf of the canonical embedding $C
\subset \PP^{g - 1}$ starts as
\[ \bigoplus_{i \geq 0} \O_{\PP^{g - 1}} (-3 - i)^{\oplus b_i} \hpil \O_{\PP^{g - 1}} (-2)^{\oplus a} \hpil
\I_{C/\PP^{g - 1}} \hpil 0 \]  
with $b_i = 0$ for $i \geq 1$ when $\Cliff(C) \geq 3$. Restricting to $C$ and dualizing we get an exact
sequence
\[ 0 \hpil N_{C / \PP^{g - 1}} \hpil \O_C (2)^{\oplus a} \hpil \bigoplus_{i \geq 0} \O_C (3
+ i)^{\oplus b_i} \]
whence an exact sequence
\[ 0 \hpil H^0(N_{C / \PP^{g - 1}} \otimes M^{-1}) \hpil H^0(2K_C - M)^{\oplus a}
\mapright{\varphi} \bigoplus_{i \geq 0} H^0((3+i)K_C - M)^{\oplus b_i}. \]
As in the proof of Lemma \ref{N2} we have that $H^0(N_{C / \PP^{g - 1}} \otimes M^{-1}) = 0$ under hypothesis
(i) and $H^0(N_{C / \PP^{g - 1}} \otimes M^{-1}) = \Ker \varphi = 0$ under hypothesis (ii). Therefore
we conclude by Proposition \ref{gaussian} (a).
\end{proof}

\noindent Using an appropriate generalization of the methods of \cite[Proof of Thm.2]{bel} we can also get
surjectivity when $h^0(2K_C - M) \geq 2$.

\begin{prop}
\label{bel}
Let $C$ be a smooth irreducible curve of genus $g \geq 4$ and let $M$ be a line bundle on $C$. Suppose there 
exists an integer $m \geq 1$ and an effective divisor $D = P_1 + \ldots + P_m$ such that
\begin{itemize}
\item[(i)] $H^1(M - 2P_i ) = 0$ for $1 \leq i \leq m$;
\item[(ii)] $h^0(D) = 1$ and $h^0(2K_C - M -D) = 0$;
\item[(iii)] $m \leq \Cliff(C) - 2$.
\end{itemize}
Then $\Phi_{M, \omega_C}$ is surjective.
\end{prop}

\begin{proof}
As is well-known we have $\Cliff(C) \leq \lfloor \frac{g - 1}{2} \rfloor$, whence $m \leq \Cliff(C) - 2 \leq g - 4$. We
start by observing that $K_C - D$ is very ample. In fact, if $K_C - D$ is not very ample, there are two points $Q_1, Q_2 \in
C$ such that 
\[ h^0(K_C - D - Q_1 - Q_2) = h^0(K_C - D) - 1 = g - 2 - m + h^0(D) = g - 1 - m \]
whence $h^1(D + Q_1+ Q_2) = g - 1 - m \geq 3$ and $h^0(D + Q_1+ Q_2) = 2$ by Riemann-Roch. Therefore
$D + Q_1+ Q_2$  contributes to the Clifford index of $C$ and we have $\Cliff(C) \leq \Cliff(D + Q_1+ Q_2) =
m$, contradicting (iii). 

\noindent Consider the embedding $C \subset \PP H^0(K_C - D) = \PP^r$, where $r = g - 1 - m$. We claim that,
in the latter embedding, $C$  has no trisecant lines. As a matter of fact if there exist three points $Q_1,
Q_2, Q_3  \in C$ such that their linear span $<Q_1, Q_2, Q_3>$ is a line, we have that
\begin{eqnarray*}
1 & = & \dim<Q_1, Q_2, Q_3> = h^0(K_C - D) - 1 - h^0(K_C - D - Q_1 - Q_2 - Q_3) = \\
& = & g - 1 - m  - h^0(K_C - D - Q_1 - Q_2 - Q_3)
\end{eqnarray*}
whence $h^1(D + Q_1+ Q_2 + Q_3) = g - 2 - m \geq 2$ and again $h^0(D + Q_1+ Q_2 + Q_3) = 2$. 
Therefore $D + Q_1+ Q_2 + Q_3$ contributes to the Clifford index of $C$ and we get $\Cliff(C) \leq 
\Cliff(D + Q_1+ Q_2 + Q_3) = m + 1$, contradicting (iii).

\noindent Note further that by (ii) and (iii) we have
\[ \deg(K_C - D) = 2g - 2 - m \geq 2g + 2 - 2h^1(K_C - D) - \Cliff(C) \]
therefore Green-Lazarsfeld's theorem \cite[Prop.2.4.2]{la1} gives that $C$ is scheme-theoretically cut out by 
quadrics in $\PP^r$. Hence we have a surjection 
\[ \O_C(2D - 2K_C)^{\oplus \alpha} \to N_{C / \PP^r}^{\ast} \to 0. \]
Setting, as in \cite{bel}, $R_L = N_{C / \PP H^0(L)}^{\ast} \otimes L$ for any very ample line bundle 
$L$, we deduce a surjection
\[ \O_C(M - K_C + D)^{\oplus \alpha} \to R_{K_C - D} \otimes M \to 0. \]
By (ii), we have that $H^1(M - K_C + D) = H^0(2K_C - M - D)^{\ast} = 0$ whence 
\[ H^1(R_{K_C - D} \otimes M) = 0. \]
Now there is an exact sequence \cite[2.7]{bel}, \cite[Proof of Thm.5]{ei}
\[ 0 \hpil R_{K_C - D} \otimes M \hpil R_{K_C} \otimes M \hpil \bigoplus\limits_{i = 1}^m \O_C(M - 2P_i ) 
\hpil 0 \]
and therefore by (i) we deduce that 
\[H^0(N_{C / \PP^{g-1}} \otimes M^{-1}) \cong H^1(N_{C / \PP^{g-1}}^{\ast} \otimes \omega_C 
\otimes M)^{\ast} \cong H^1(R_{K_C} \otimes M)^{\ast} = 0. \] 
Hence we get the surjectivity of $\Phi_{M, \omega_C}$ by Proposition \ref{gaussian} (a).
\end{proof}

\noindent We will often use the above result in the following simplified version.

\begin{cor}
\label{bel2}
Let $C$ be a smooth irreducible curve of genus $g \geq 4$ and let $M$ be a line bundle on $C$ such that 
$H^1(M) = 0$ and $\deg(M) \geq g + 1$. Suppose that
\[ h^0(2K_C - M) \leq \Cliff(C) - 2. \]
Then $\Phi_{M, \omega_C}$ is surjective.
\end{cor}

\begin{proof}
Let $m = \Cliff(C) - 2$. Then $m \geq 0$ by hypothesis and when $m = 0$ the surjectivity of $\Phi_{M,
\omega_C}$ holds by Proposition \ref{Cliff}. When $m \geq 1$ choose general points $P_1, \ldots, P_m$ of $C$ 
and apply Proposition \ref{bel}. 
\end{proof}

\begin{cor} \cite[Cor.1.7]{te}
Let $C$ be a smooth irreducible curve of genus $g \geq 5$ nontrigonal and not isomorphic to a plane
quintic. Let $M$ be a line bundle on $C$. 

\noindent Then the Gaussian map $\Phi_{M, \omega_C}$ is surjective if $\deg M \geq 4g - 4$ and $M \neq 2K_C$ if
equality holds.
\end{cor}

\begin{proof}
Immediate consequence of Corollary \ref{bel2} or of Proposition \ref{Cliff}.
\end{proof}

\subsection{Gaussian maps on tetragonal curves}
\label{tetra}

In this subsection we improve Tendian's \cite{te} results about Gaussian maps on tetragonal curves.
Moreover note that, even though the statement in \cite[Thm.2.10]{te} is almost correct, the proof certainly 
contains a gap (see Remark \ref{ten}).

\noindent We start with some generalities on tetragonal curves following again \cite[6.2]{sch1}. 

\begin{defn-not}
\label{b2}
Let $C$ be a smooth irreducible tetragonal curve of genus $g \geq 6$ not isomorphic to a plane quintic. Let 
$A$ be a $g^1_4$ on $C$ and let $V_A \subset \PP^{g-1} = \PP H^0(\omega_C)$ be the rational normal scroll
spanned by the divisors in $|A|$, $H_A$ the hyperplane bundle and $R_A$ a ruling of $V_A$. Let $\E_A$ be the 
rank 3 vector bundle on $\PP^1$ so that $V_A$ is the image of $\PP \E_A$ under the morphism given by $|\O_{\PP
\E_{A}}(1)|$. Let $\widetilde{H}_A$ and $\widetilde{R}_A$ be the pull-backs, under this morphism, of $H_A$ and
$R_A$ respectively. Then there are two integers $b_{1,A}, b_{2,A}$ such that $b_{1,A} \geq b_{2,A} \geq 0$,
$b_{1,A} + b_{2,A} = g - 5$ and there are two surfaces $\widetilde{Y}_A \sim 2 \widetilde{H}_A - b_{1,A}
\widetilde{R}_A$, $\widetilde{Z}_A \sim 2 \widetilde{H}_A -  b_{2,A} \widetilde{R}_A$ such that, if $Y_A, Z_A$
are their images in $\PP^{g-1}$ then $C = Y_A \cap Z_A$. We also define 
\[ b_2(C) = \min\{b_{2,A}, A \ \mbox{a} \ g^1_4 \ \mbox{on} \ C \}. \]
\end{defn-not}

\noindent We have

\begin{lemma}
\label{N2perY}
The surface $Y_A \subset \PP^{g-1}$ has degree $g - 1 + b_{2, A}$ and satisfies property $N_2$.
\end{lemma}

\begin{proof}
We set for simplicity $Y = Y_A$, $\widetilde{Y} = \widetilde{Y}_A$, $V = V_A$, $\E = \E_A$, $H = H_A$, $R = R_A$, 
$\widetilde{H} = \widetilde{H}_A$, $\widetilde{R} = \widetilde{R}_A$, $b_i = b_{i,A}, i = 1, 2$. Note that
$\widetilde{H}^3 = \deg V = g - 3$, $\widetilde{R}^2 = 0$ and $\widetilde{H}^2.\widetilde{R} = 1$. Let
$\widetilde{X} \in |\O_{\widetilde{Y}}(\widetilde{H})|$ be a general curve. Since
$|\O_{\widetilde{Y}}(\widetilde{H})|$ is not composed with a pencil we have that $\widetilde{X}$ is irreducible.
Moreover $\widetilde{X}$ is smooth outside $\widetilde{H} \cap \Sing(\widetilde{Y})$, whence $\widetilde{X}$ is
also reduced. 

\noindent Let $\L = \O_{\widetilde{X}}(\widetilde{H}), X = Y \cap H$, so that $\varphi_{\L}(\widetilde{X}) = X$.
We will first prove that $X$ satisfies property $N_2$. 

\noindent To this end by \cite[Thm.A]{bf} it is enough to show that 
\begin{equation}
\label{ineq}
\deg X \geq 2p_a(X) + 3.
\end{equation}
Taking intersections in $\PP \E$ we have
\begin{equation}
\label{degree}
\deg X = \deg \L = \widetilde{H}^2 \cdot \widetilde{Y} = \widetilde{H}^2 \cdot (2\widetilde{H} - b_1 \widetilde{R})
= 2g - 6 - b_1 = g - 1 + b_2.
\end{equation}
On the other hand, using the cohomology of the scroll, we get
\begin{eqnarray*}
p_a(X) & = & 1 - \chi(\O_X) = 1 - \chi(\O_Y) + \chi(\O_Y(-1)) = \\ 
& = & 1 - \chi(\O_V) + \chi(\O_V(-2 H + b_1 R)) + \chi(\O_V(-1)) - \chi(\O_V(-3 H + b_1 R)) = \\
& = & g - 4 - b_1.
\end{eqnarray*}

\noindent Now $2p_a(X) + 3 = 2g - 5 - 2b_1 \leq 2g - 6 - b_1$ if and only if $b_1 \geq 1$. The latter holds
because $b_1 \geq b_2 \geq 0$ and $g \geq 6$. Therefore (\ref{ineq}) is proved.

\noindent Again using the cohomology of the scroll it is easy to prove that $H^1(\I_{Y / \PP^{g-1}}(j)) = 0$ for
every $j \in \ZZ$ and that $H^1(\O_Y(j)) = 0$ for every $j \geq 0$. Applying \cite[Thm.2.a.15 and Thm.3.b.7]{gr}
(that hold for any scheme) we deduce that $Y$ satisfies property $N_2$ since $Y \cap H$ does.
\end{proof}

\begin{rem}
\label{ten} 
{\rm In Tendian's paper it is assumed that a general hyperplane section $Y \cap H$ is smooth, but 
in fact it can be singular \cite[6.5]{sch1} when the $g^1_4$ exhibits $C$ as a double cover of an elliptic or
hyperelliptic curve.}
\end{rem} 

\begin{prop}
\label{tetragonal}
Let $C$ be a smooth irreducible tetragonal curve of genus $g \geq 6$ not isomorphic to a plane quintic. Let $A$
be a $g^1_4$, set $b_2 = b_{2, A}$ and let $M$ be a line bundle on $C$. We have
\begin{itemize}
\item[(i)] If $h^0(2K_C - M) \leq 1$ and $h^0(2K_C - M - b_2 A) = 0$, then $\Phi_{M, \omega_C}$ is surjective;
\item[(ii)] If $H^1(M) = 0$ and $\mu_{M, \omega_C}$ is surjective, then 
$\cork \Phi_{M, \omega_C} \geq h^0(2K_C - M - b_2 A)$, with equality holding if $h^0(2K_C - M) \leq 1$.
\end{itemize}
\end{prop}

\begin{proof}
Let $Y$ be the surface arising in the scroll defined by $A$ and set, as in Lemma \ref{N2}, $\F_{C,
Y} = \Shom_{\O_{\PP^{g-1}}}(\I_{Y / \PP^{g-1}}, \O_C)$. Applying the left exact functor 
$\Shom_{\O_{\PP^{g-1}}}(- , \O_C)$ to the exact sequence
\[ 0 \hpil \I_{Y / \PP^{g-1}} \hpil \I_{C / \PP^{g-1}} \hpil \I_{C / Y} \hpil 0 \]
we get an exact sequence
\begin{equation}
\label{normal} 
0 \hpil N_{C / Y} \otimes M^{-1} \hpil N_{C / \PP^{g-1}} \otimes M^{-1} \hpil \F_{C, Y} \otimes M^{-1}.
\end{equation}
Observe that $h^0(N_{C / Y} \otimes M^{-1} ) = h^0(2K_C - M - b_2 A)$. Now if $h^0(2K_C - M - b_2 A) = 0$,
from (\ref{normal}) and Proposition \ref{gaussian} (a), we see that to prove (i) we just need to show that
\begin{equation}
\label{restrnorm} 
\mbox{if} \  h^0(2K_C - M) \leq 1 \ \mbox{then} \ H^0(\F_{C, Y} \otimes M^{-1}) = 0.
\end{equation}
On the other hand, under the hypotheses in (ii), we have that $\cork \Phi_{M, \omega_C} = 
h^0(N_{C / \PP^{g-1}} \otimes M^{-1})$ by Proposition \ref{gaussian} (b). Now from
(\ref{normal}) we get that $h^0(N_{C / \PP^{g-1}} \otimes M^{-1}) \geq h^0(2K_C - M - b_2 A)$ and to prove
equality we need again to prove (\ref{restrnorm}).

\noindent To conclude we just note that (\ref{restrnorm}) holds by Lemmas \ref{N2} and \ref{N2perY}.
\end{proof}

\section{Linear series on quadric sections of surfaces of degree $g - 1$ in $\PP^{g - 1}$}

\noindent In this section we will use some well-known vector bundle methods (\cite{la1}, \cite{ty}) to study 
linear series on curves of genus $g$ that are, in their canonical embedding, a quadric section of a surface of
degree $g - 1$ in $\PP^{g - 1}$. We recall that when the surface is a smooth Del Pezzo the gonality and Clifford index of
such curves are known by \cite{pa}, \cite{kn1}. Most of the results we prove are probably known, at least in the
smooth case, but we include them anyway for completeness' sake.

\begin{lemma}
\label{bogreider}
Let $X$ be a smooth surface with $- K_X \geq 0$. Let $C \subset X$ be a smooth irreducible curve of genus $g$ 
and let $A$ be a base-point free $g^1_k$ on $C$. Suppose that $2g - 2 - K_X.C - 4k \geq \max\{0, 3 - 4 \chi(\O_X)
\}$ and, if $h^1(\O_X) \geq 1$, that $h^0(N_{C/X} \otimes A^{-1}) \geq 2 h^1(\O_X) + 1$. Then there exist two line
bundles $L, M$ on $X$ and a zero-dimensional subscheme $Z \subset X$ such that the following hold:
\begin{itemize}
\item[(i)] $C \sim M + L$;
\item[(ii)] $k = M.L + \length(Z) \geq M.L \geq L^2 \geq 0$;
\item[(iii)] there exists an effective divisor $D$ on $C$ of degree $M.L + L^2 - k \geq 0$ such that 

\noindent $A \cong L_{|C}(- D)$;
\item[(iv)] if $L^2 = 0$ then $M.L = k$ and $A \cong L_{|C}$;
\item[(v)] $L$ is base-component free and nontrivial;
\item[(vi)] if $C \sim -2K_X$ then $3L^2 + M.L \in 4 \ZZ$.
\end{itemize}
\end{lemma}
\begin{proof}
Let $\F = \Ker \{H^0(A) \otimes \O_X \to A \}$ and $\E = \F^{\ast}$. As is well-known (\cite{la1}) $\E$ is a 
rank two vector bundle sitting in an exact sequence
\begin{equation}
\label{laz}
0 \hpil H^0(A)^{\ast} \otimes \O_X \hpil \E \hpil N_{C/X} \otimes A^{-1} \hpil 0
\end{equation}
and moreover $c_1(\E) = C$ and $c_2(\E) = k$, so that $\Delta(\E) := c_1(\E)^2 - 4c_2(\E) = C^2 - 4k =
2g - 2 - K_X.C - 4k \geq 0$. Let $H$ be an ample line bundle on $X$ and suppose that $\E$ is $H$-stable. Then
$h^0(\E \otimes \E^{\ast}) = 1$ by \cite[Cor.4.8]{fr} and $h^2(\E \otimes \E^{\ast}) = h^0(\E \otimes
\E^{\ast}(K_X)) \leq  h^0(\E \otimes \E^{\ast}) = 1$, therefore $2 \geq h^0(\E \otimes \E^{\ast}) + h^2(\E
\otimes \E^{\ast}) = h^1(\E \otimes \E^{\ast}) + \chi(\E \otimes \E^{\ast}) \geq 4 \chi(\O_X) + \Delta(\E) \geq
3$, a contradiction. Hence $\E$ is not $H$-stable and if $M$ is the maximal destabilizing subbundle we have an
exact sequence
\begin{equation}
\label{dest}
0 \hpil M \hpil \E \hpil \I_{Z/X} \otimes L \hpil 0
\end{equation}
where $L$ is another line bundle on $X$ and $Z$ is a zero-dimensional subscheme of $X$. Computing Chern
classes in (\ref{dest}) we get (i) and the equality in (ii). Since the destabilizing condition reads 
$(M - L).H \geq 0$ and since $(M - L)^2 = \Delta(\E) + 4\length(Z) \geq 0$, we see that $M - L$ belongs to
the closure of the positive cone of $X$. We now claim that $\E$ is globally generated off a finite set. In fact 
if $h^1(\O_X) \geq 1$ we have by hypothesis that $h^0(N_{C/X} \otimes A^{-1}) \geq 2 h^1(\O_X) + 1$ and the claim
follows by (\ref{laz}) since the map $\psi: H^0(\E) \to H^0(N_{C/X} \otimes A^{-1})$ is nonzero. On the other
hand if $h^1(\O_X) = 0$ we have that $\psi$ is surjective, whence, again by (\ref{laz}), we just need to prove
that $h^0(N_{C/X} \otimes A^{-1}) \geq 1$. Since $g \geq 2k + 1 + \frac{1}{2} K_X.C$ we get $\deg (N_{C/X} \otimes
A^{-1}) = 2g - 2 - K_X.C - k \geq g$. Therefore $h^0(N_{C/X} \otimes A^{-1}) \geq 1$ by Riemann-Roch and the
claim is proved.

\noindent Since $\E$ is globally generated off a finite set then so is $L$. It follows that $L \geq 0$, $L$ is 
base-component free and $L^2 \geq 0$. Now the signature theorem \cite[VIII.1]{bpv} implies that $(M - L).L \geq
0$ thus proving (ii). To see (iii) and (iv) note that if $M.L > 0$ then the nefness of $L$ implies that 
$H^0(- M) = 0$.  On the other hand if $M.L = 0$ then $L^2 = C.L = 0$ whence $L \eqv 0$ by the Hodge index theorem
and therefore $C \eqv M$. Then $M.H = C.H > 0$ whence again $H^0(-M) = 0$. Twisting (\ref{laz}) and  (\ref{dest})
by $-M$ we deduce that $h^0(L_{|C} \otimes A^{-1}) \geq h^0(\E(-M)) \geq 1$. This proves (iii) and also (v).
Moreover it gives $\deg(L_{|C} \otimes A^{-1}) \geq 0$, whence, if $L^2 = 0$, we get that $M.L \geq k$. By (ii)
it follows that $M.L = k$ and therefore $\deg(L_{|C} \otimes A^{-1}) = 0$, whence $L_{|C} \cong A$. This proves
(iv). 

\noindent Finally suppose that $C \sim -2K_X$. We have $\chi(L) = \chi(\O_X) + \frac{1}{2} L.(L - K_X)$ whence 
$2 L.(L - K_X)$ is divisible by $4$. But $2 L.(L - K_X) = 2L^2 + L.C = 3L^2 + M.L$, giving (vi).
\end{proof}

\noindent We now analyze linear series on curves on surfaces of degree $r$ in $\PP^r$. We will use the following

\begin{defn-not}
\label{delp}
For $1 \leq n \leq 9$ we denote by $\Sigma_n$ the blow-up of $\PP^2$ at $n$ possibly infinitely near
points, by $\widetilde{H}$ the strict transform of a line and by $G_i$ the total inverse image of the blown-up 
points. Let $Q \subset \PP^3$ be a quadric cone with vertex $V$. We denote by $Bl_{V}Q$ the blow-up of $Q$ along
$V$ and by $\widetilde{H}$ the strict transform of a plane. Let $C_n \subset \PP^n$ be the cone over a smooth
elliptic curve in $\PP^{n-1}$ and let $V$ be the vertex. We denote by $Bl_{V} C_n$ the blow-up of $C_n$ along
$V$, by $C_0$ the inverse image of $V$ and by $f$ the numerical class of a fiber.
\end{defn-not}

\begin{rem}
\label{delp2}
{\rm We recall that by \cite[Thm.8]{na} a linearly normal integral surface $Y \subset \PP^r$ of degree $r$ is
either the anticanonical image of $\Sigma_{9-r}$ or $C_r$ or the 2-Veronese embedding in $\PP^8$ of an
irreducible quadric in $\PP^3$ or the 3-Veronese embedding in $\PP^9$ of $\PP^2$.}
\end{rem} 

\begin{prop}
\label{g1kondelp}
Let $X$ be a surface among $\Sigma_n$, $Bl_{V}Q$ or $Bl_{V} C_n$ as in Definition {\rm \ref{delp}} and let $C$ be a
smooth irreducible curve such that, if $X = \Sigma_n$ or $Bl_{V}Q$ then $C \sim -2K_X$, while if $X = Bl_{V}
C_n$ then $C \eqv -2K_X - 2C_0$. We have:
\begin{itemize}
\item[(a)] if $X = \Sigma_1$ then $C$ has no complete base-point free $g^1_6$;
\item[(b)] if $X = \Sigma_2$ then every complete base-point free $g^1_4$ on $C$ is $(\widetilde{H} -
G_i)_{|C}$, $i = 1, 2$;
\item[(c)] if $X = \Sigma_2$ then every complete base-point free $g^1_6$ on $C$ is $(2\widetilde{H} -
G_1 - G_2)_{|C} - P_1 - P_2$, where $P_1, P_2$ are two points of $C$;
\item[(d)] if $X = \Sigma_3$ then every complete base-point free $g^1_4$ on $C$ is $(\widetilde{H} -
G_i)_{|C}$, $i = 1, 2, 3$;
\item[(e)] if $X = \Sigma_3$ then every complete base-point free $g^1_5$ on $C$ is either
$\widetilde{H}_{|C} - P$ or $(2\widetilde{H} - G_1 - G_2 - G_3)_{|C} - P$, for some point $P \in C$;
\item[(f)] if $X = \Sigma_3$ and $A$ is a complete base-point free $g^1_6$ on $C$ then either
$A \cong (2\widetilde{H} - G_i - G_j)_{|C} - P_1 - P_2$, for $1 \leq i < j \leq 3$ and $P_1, P_2$ are two
points of $C$ or $(-K_X)_{|C} - A$ is another complete base-point free $g^1_6$ on $C$ different from 
$(2\widetilde{H} - G_i - G_j)_{|C} - P_1 - P_2$;
\item[(g)] if $X = Bl_{V} C_6$ then $C$ has no complete base-point free $g^1_5$ and every complete 
base-point free $g^1_4$ on $C$ is $(f_1 + f_2)_{|C}$, where $f_1, f_2$ are two fibers;
\item[(h)] if $X = Bl_{V}Q$ then $C$ has a unique complete base-point free $g^1_4$, namely $f_{|C}$, where
$f$ is the pull-back of a line of the cone $Q$;
\item[(i)] if $X = Bl_{V}Q$ then every complete base-point free $g^1_6$ on $C$ is $\widetilde{H}_{|C} - P_1
- P_2$, where $P_1, P_2$ are two points of $C$;
\item[(j)] if $X = Bl_{V}Q$ then there is no effective divisor $Z \subset C$ such that $f_{|C} + Z$ is a
complete base-point free $g^2_8$ on $C$.
\end{itemize}
\end{prop}

\begin{proof} We record, for later use, the following fact on $X = \Sigma_n$. Let $\L$ be a nef line bundle
on $X$ with $\L \sim a \widetilde{H} - \sum\limits_{i = 1}^n b_i G_i$. Then 
\begin{equation} 
\label{cs1}
a = \L.\widetilde{H} \geq 0, \ \ b_i = \L.G_i \geq 0, \ \ \L^2 = a^2 - \sum\limits_{i = 1}^n b_i^2, \ \ \L.(- K_X) 
= 3a - \sum\limits_{i = 1}^n b_i 
\end{equation} 
and the Cauchy-Schwartz inequality $(\sum\limits_{i = 1}^n b_i)^2 \leq n \sum\limits_{i = 1}^n b_i^2$ implies
that
\begin{equation}
\label{cs2}
(3a + \L.K_X)^2 \leq n (a^2 - \L^2).
\end{equation} 

\noindent We will now apply Lemma \ref{bogreider} to a base-point free $g^1_k$ indicated in (a)-(i) and we
will set $z = \length(Z)$. 

\noindent {\bf (a)} We have $K_X^2 = 8$ whence $C^2 = 32, k = 6$ and from (ii) of Lemma \ref{bogreider} we deduce
that $6 = M.L + z \geq M.L \geq L^2 \geq 0$. Now if $3 \leq L^2 \leq 6$ we have a contradiction by the Hodge
index theorem applied to $C$ and $L$. The same theorem implies, for $L^2 = 2$, that $C \eqv 4L$. But $C
\sim 6\widetilde{H} - 2G_1$ whence the contradiction $4L.\widetilde{H} = C.\widetilde{H} = 6$. If $L^2 = 1$
write $L \sim a\widetilde{H} - b_1G_1$. Then $a^2 = b_1^2 + 1$ therefore $a = 1, b_1 = 0$ and $L \sim
\widetilde{H}$. Then $\deg(L_{|C} \otimes A^{-1}) = \widetilde{H}.C - 6 = 0$, whence $A \cong
\widetilde{H}_{|C}$ by (iii) of Lemma \ref{bogreider}. Therefore we have the contradiction $h^0(A) = 3$. If
$L^2 = 0$ by (iv) of Lemma \ref{bogreider} we have that $M.L = 6$ whence $3L^2 + M.L = 6$, contradicting (vi)
of Lemma \ref{bogreider}. This proves (a).

\noindent {\bf (b)} We have $K_X^2 = 7, C^2 = 28$ and $k = 4$. By (ii) of Lemma \ref{bogreider} and the Hodge 
index theorem applied to $C$ and $L$ we see that we are left with the case $L^2 = 0$ whence $A \cong L_{|C}$. By
(\ref{cs1}), (\ref{cs2}) we deduce that $L \sim \widetilde{H} - G_i$ for $i = 1, 2$.  This proves (b).

\noindent {\bf (c)} We have $K_X^2 = 7$ whence $C^2 = 28$ and $k = 6$. From (ii) of Lemma \ref{bogreider} and the
Hodge index theorem applied to $C$ and $L$ we get $0 \leq L^2 \leq 2$. The same theorem implies, for $L^2 =
2$, that $z = 0, M.L = 6$. By (iii) of Lemma \ref{bogreider} we have that there are two points $P_1, P_2 \in 
C$ such that $A \cong L_{|C} - P_1 - P_2$. By (\ref{cs1}), (\ref{cs2}) we deduce that $L \sim 2\widetilde{H} -
G_1 - G_2$. If $L^2 = 1$ again by (ii) of Lemma \ref{bogreider} and the Hodge index theorem applied to $C$ and
$L$ we get that $0 \leq z \leq 1$ and $5 \leq M.L \leq 6$. By (vi) of Lemma \ref{bogreider} we have that $M.L = 5$
whence $\deg(L_{|C} \otimes A^{-1}) = 0$, so that $A \cong L_{|C}$ by (iii) of Lemma \ref{bogreider}. By
(\ref{cs1}), (\ref{cs2}) we deduce that $L \sim \widetilde{H}$, giving the contradiction $h^0(A) = 3$. If $L^2 =
0$ we have that $M.L = 6$ by (iv) of Lemma \ref{bogreider} contradicting (vi) of Lemma \ref{bogreider}. This
proves (c).

\noindent {\bf (d)} We have $K_X^2 = 6$ whence $C^2 = 24$ and $k = 4$. From (ii) of Lemma \ref{bogreider} and the
Hodge index theorem applied to $C$ and $L$ we get $0 \leq L^2 \leq 1$. The same theorem implies, for $L^2 =
1$, that $z = 0, M.L = 4$, contradicting (vi) of Lemma \ref{bogreider}. Therefore $L^2 = 0$ and (iv) of 
Lemma \ref{bogreider} implies that $M.L = 4$ and $A \cong L_{|C}$. By (\ref{cs1}), (\ref{cs2}) we deduce that $L
\sim \widetilde{H} - G_i$. This proves (d).

\noindent {\bf (e)} We have $K_X^2 = 6$ whence $C^2 = 24$ and $k = 5$. From (ii) of Lemma \ref{bogreider} and the
Hodge index theorem applied to $C$ and $L$ we get $L^2 \leq 2$ with equality only when $z = 0, M.L = 5$,
contradicting (vi) of Lemma \ref{bogreider}. When $L^2 = 1$, the same theorem together with (vi) of
Lemma \ref{bogreider} implies that $z = 0, M.L = 5$, whence $A \cong L_{|C} - P$ by (iii) of Lemma
\ref{bogreider}. By (\ref{cs1}), (\ref{cs2}) we deduce that either $L \sim \widetilde{H}$ or $L \sim
2\widetilde{H} - G_1 - G_2 - G_3$. If $L^2 = 0$ then (iv) of Lemma \ref{bogreider} implies that $M.L = 5$,
contradicting (vi) of Lemma \ref{bogreider}. This proves (e).

\noindent {\bf (f)} We have $K_X^2 = 6$ whence $C^2 = 24$ and $k = 6$. From (ii) of Lemma \ref{bogreider} and the
Hodge index theorem applied to $C$ and $L$ we see, for $3 \leq L^2 \leq 5$, that $z = 0, M.L = 6$, 
contradicting (vi) of Lemma \ref{bogreider}. If $L^2 = 2$ by the Hodge index theorem and (vi) of Lemma
\ref{bogreider} we have that $z = 0, M.L = 6$. By (\ref{cs1}), (\ref{cs2}) we deduce that $L \sim 2\widetilde{H}
- G_i - G_j$ for $i \neq j$ and by (iii) of Lemma \ref{bogreider} we have that there are two points $P_1, P_2
\in C$ such that $A \cong L_{|C} - P_1 - P_2$. If $L^2 = 1$ by the Hodge index theorem and (vi) of Lemma
\ref{bogreider} we have that $z = 1, M.L = 5$. By (\ref{cs1}), (\ref{cs2}) we deduce that either $L \sim
\widetilde{H}$ or $L \sim 2\widetilde{H} - G_1 - G_2 - G_3$. By (iii) of Lemma \ref{bogreider} we have that $A
\cong L_{|C}$,  giving the contradiction $h^0(A) = 3$. If $L^2 = 0$ by (iv) of Lemma \ref{bogreider} we have
that $M.L = 6$ contradicting (vi) of Lemma \ref{bogreider}. Finally when $L^2 = 6$ the Hodge index theorem
applied to $C$ and $L$ implies that $C \eqv 2L$ and $z = 0$. Therefore $L \sim M \sim - K_X$ whence the exact
sequence (\ref{dest}) splits since $\Ext^1(\O_X(- K_X), \O_X(- K_X)) = 0$ and we get $\E \cong \O_X(-
K_X)^{\oplus 2}$. Therefore $\E(K_X)$ is globally generated and so is $(- K_X)_{|C} \otimes A^{-1}$ by
(\ref{laz}). Moreover again by (\ref{laz}) we get that $(- K_X)_{|C} \otimes A^{-1}$ is a $g^1_6$. Also such
a $g^1_6$ cannot coincide with the other type $(2\widetilde{H} - G_i - G_j)_{|C} - P_1 - P_2$, for otherwise we
would have that $(- K_X)_{|C} \otimes A^{-1} \sim (2\widetilde{H} - G_i - G_j)_{|C} - P_1 - P_2$, whence $A
\cong (\widetilde{H} - G_k)_{|C} + P_1 + P_2$ would have two base points. This proves (f).
 
\noindent {\bf (g)} We have that $X \cong \PP(\O_E \oplus \O_E(-1))$ where $E \subset \PP^5$ is a smooth elliptic
normal curve. Let $C_0$ be a section and $f$ be a fiber so that $C_0^2 = - 6$ and the intersection form is even.
Moreover $C \eqv 2C_0 + 12 f$, $C^2 = 24$ and $k = 4, 5$.  From (ii) of Lemma \ref{bogreider} and the Hodge index 
theorem applied to $C$ and $L$ we deduce, if $L^2 \geq 2$, that $k = 5$, $L^2 = 2$, $z = 0$ and $M.L = 5$. On the
other hand if $L^2 = 0$ we have that $M.L = k$ and $A \cong L_{|C}$ by (iv) of Lemma \ref{bogreider}. Let $L \eqv
aC_0 + bf$ so that $M \eqv (2 - a) C_0 + (12 - b)f$ and $L^2 = 2a(b - 3a)$. Moreover, by (v) of Lemma
\ref{bogreider} we have $a = f.L \geq 0$. Now if $L^2 = 2$ we get $a = 1$, $b = 4$ giving the contradiction $M.L
= 6$. Therefore $L^2 = 0$ whence either $a = 0$ or $b = 3a$. In the second case we get $k = M.L = 6a$, a
contradiction. Therefore $a = 0$ and $k = M.L = 2b$, that is $k = 4$, $b = 2$ and $L \eqv
2f$ as desired. This proves (g).

\noindent {\bf (h)} We have that $X \cong \PP (\O_{\PP^1} \oplus \O_{\PP^1}(-2))$. Let $C_0$ be a section and
$f$ be a fiber so that $C_0^2 = -2$ and the intersection form is even. Moreover $C \sim 4C_0 + 8 f$, $C^2 =
32$ and $k = 4$.  From (ii) of Lemma \ref{bogreider} and the Hodge index theorem applied to $C$ and $L$ we
have a contradiction if $L^2 \geq 2$. Hence $L^2 = 0, M.L = 4$ and $A \cong L_{|C}$ by (iv) of Lemma
\ref{bogreider}. Then we get that either $L \sim f$ or $L \sim C_0 + f$. Since $C_0.C = 0$, this proves (h).

\noindent {\bf (i)} We retain the notation used in (h) except that now $k = 6$. From (ii) of Lemma
\ref{bogreider} and  the Hodge index theorem applied to $C$ and $L$ we deduce, if $L^2 \geq 2$, that $L^2 = 2$,
$z = 0, M.L = 6$ and $C \eqv 4L$, whence $L \sim C_0 + 2 f \sim \widetilde{H}$. By (iii) of Lemma \ref{bogreider}
we have that there are two points $P_1, P_2 \in C$ such that $A \cong \widetilde{H}_{|C} - P_1 - P_2$. When $L^2
= 0$ we get $M.L = 6$ by (iv) of Lemma \ref{bogreider}, contradicting (vi) of Lemma \ref{bogreider}. This proves
(i).

\noindent {\bf (j)} Again we use the notation in (i). Suppose there is an effective divisor $Z \subset C$ such
that $f_{|C} + Z$ is a complete base-point free $g^2_8$ on $C$. By Riemann-Roch we get that 
\[ h^0((2C_0 + 3f)_{|C} - Z) = h^0(K_C - f_{|C} - Z) = 3 \] 
and the exact sequence
\[ 0 \hpil \O_X(- 2C_0 - 5f) \hpil \I_{Z/X}(2C_0 + 3f) \hpil \I_{Z/C}(2C_0 + 3f) \hpil 0 \]
gives that also $h^0(\I_{Z/X}(2C_0 + 3f)) = 3$, whence, since $h^0(2C_0 + 3f) = 6$, that $Z$ does not impose
independent conditions to $|2C_0 + 3f|$. Now let $Z' \subset Z$ be an effective divisor of degree $3$ and set
$Z' + P = Z$. By the exact sequence
\[ 0 \hpil \O_X(- 2C_0 - 5f) \hpil \I_{Z'/X}(2C_0 + 3f) \hpil \I_{Z'/C}(2C_0 + 3f) \hpil 0 \]
and Riemann-Roch we have 
\begin{eqnarray*}
h^0(\I_{Z'/X}(2C_0 + 3f)) & = & h^0(\I_{Z'/C}(2C_0 + 3f)) = h^1(f_{|C} + Z - P) = \\ 
& = & h^0(f_{|C} + Z - P) + 1 = 3.
\end{eqnarray*}
Therefore $Z$ is in special position with respect to $2C_0 + 3f \sim \L + K_X$, where $\L \sim 4C_0 + 7f$. By
\cite{re}, \cite{gh}, \cite{ca}, \cite{la2} there is a rank $2$ vector bundle $\E$ on $X$ sitting in an exact
sequence
\begin{equation}
\label{rei}
0 \hpil \O_X \hpil \E \hpil \I_{Z/X} \otimes \L \hpil 0
\end{equation}
with $c_1(\E) = \L$ and $c_2(\E) = 4$ so that $\Delta(\E) = \L^2 - 16 = 8 > 0$. Therefore $\E$ is Bogomolov 
unstable and (\cite{bo}, \cite{re}) there are two line bundles $A, B$ on $X$ and a zero-dimensional subscheme 
$W \subset X$ sitting in an exact sequence
\begin{equation}
\label{dest2}
0 \hpil A \hpil \E \hpil \I_{W/X} \otimes B \hpil 0.
\end{equation}
Moreover $\L \sim A + B$, $A.B + \length(W) = 4$, $(A - B)^2 = 8 + 4 \length(W)$ and $A - B$ lies in the
positive cone of $X$.

\noindent We record for later use two extra properties of $A$ and $B$. 

\noindent For every nef line bundle $M$ such that $M^2 \geq 0$ we have:
\begin{eqnarray}
\label{eq1} & (A - B).M \geq 0; \\ 
\label{eq2} & A.M \geq 0.
\end{eqnarray}
To prove (\ref{eq1}) and (\ref{eq2}) let $M$ be a nef line bundle such that $M^2 \geq 0$. Then $M.H \geq
0$ for every ample $H$, whence $M$ lies in the closure of the positive cone of $X$, therefore $(A - B).M \geq
0$ by \cite[VIII.1]{bpv}. Now if $A.M < 0$ then also $B.M < 0$ by (\ref{eq1}), whence $h^0(A) = h^0(B) = 0$, as
$M$ is nef. But this and (\ref{dest2}) give $h^0(\E) = 0$, contradicting (\ref{rei}).

\noindent Now $(A - B)^2 \geq 8$ and $(A + B)^2 = \L^2 = 24$ therefore
\begin{equation}
\label{ab}
A^2 + B^2 \geq 16.
\end{equation}
Moreover $\L$ lies in the positive cone of $X$, whence, by \cite[VIII.1]{bpv}, $(A - B).\L > 0$, that is
\begin{equation}
\label{ab2}
A^2 > B^2.
\end{equation}
Now if $A^2 \leq 8$ we deduce by (\ref{ab2}) that $B^2 \leq 6$, contradicting (\ref{ab}). Therefore
\begin{equation}
\label{ab3}
A^2 \geq 10.
\end{equation}
Suppose that $A \sim aC_0 + a_1f$ so that $B \sim (4 - a)C_0 + (7 - a_1)f$. Intersecting $A$ with the nef
divisors $f, C_0 + 2f$ and using (\ref{eq2}), we see that $a \geq 0, a_1 \geq 0$, whence $A \geq 0$ and in fact 
$A > 0$ by (\ref{ab3}). Also $a > 0$, for otherwise $A^2 = 0$. Now the exact sequences (\ref{rei}) and
(\ref{dest2}) twisted by $- A$ give 
\begin{equation}
\label{bz}
h^0(\I_{Z/X}(B)) \geq h^0(\E(- A)) \geq 1
\end{equation}
whence also $B > 0$. The nefness of $C_0 + 2f$ then implies $7 - a_1 = B.(C_0 + 2f) \geq 0$, whence $a_1 \leq
7$, while the nefness of $f$ implies that $4 - a = B.f \geq 0$, whence $a \leq 4$. By (\ref{eq1}) with $M =
C_0 + 2f$ we get $2 a_1 - 7 = (A - B).(C_0 + 2f) \geq 0$, whence $a_1 \geq 4$. Finally by (\ref{ab3}) we have
that $a(a_1 - a) \geq 5$. Therefore we have proved that
\begin{equation}
\label{ab4}
1 \leq a \leq 4, \ 4 \leq a_1 \leq 7, \ a(a_1 - a) \geq 5.
\end{equation}
If $a = 1, 2$ we get that $A^2 + B^2 \leq 12$, contradicting (\ref{ab}). Recall now that $C_0 \cap C =
\emptyset$ since $C_0.C = 0$. When $a = 3$ we have $A^2 + B^2 = 4a_1 - 6$ whence $a_1 = 6, 7$ by (\ref{ab}). 
When $a_1 = 7$ we have $B \sim C_0$, whence $B = C_0$. By (\ref{bz}) we deduce the contradiction $Z \subset C_0
\cap C = \emptyset$. When $a_1 = 6$ we have $B \sim C_0 + f$, whence $B = C_0 \cup F$ for some ruling $F$. As
above we have that $Z \cap C_0 = \emptyset$, whence $Z \subset F \cap C$. Since $F.C = 4$ we have that $Z = F
\cap C$, whence $Z \sim f_{|C}$ and therefore $f_{|C} + Z \sim 2f_{|C}$ is a complete base-point free $g^2_8$ on
$C$. This is of course a contradiction since on $X$ we have that $2f_{|C}$ is a complete base-point free $g^3_8$
on $C$. Finally when $a = 4$ we have $B \sim (7 - a_1)f$ whence $a_1 \leq 6$ as $B > 0$. By (\ref{ab4}) we get
$a_1 = 6$ whence $B \sim f$, therefore again $B = F$ for some ruling $F$. Hence $Z \subset F \cap C$, giving the
same contradiction above. This proves (j).
\end{proof}

\begin{rem}
{\rm Let $C$ be a smooth tetragonal curve of genus $7$ such that $\dim W^1_4(C) = 0$ and $\dim W^1_5(C) = 1$ (as
in the case $C \sim - 2 K_X$ on $X = \Sigma_3$). By \cite{acgh} $W^1_6(C)$ has an irreducible component of
dimension at least $3$ and whose general element $A$ is a complete $g^1_6$ on $C$. Moreover $A$ is base-point
free since $\dim W^1_4(C) = 0$ and $\dim W^1_5(C) = 1$. Also the same holds for $K_C - A$ thus proving that, for
these curves, there is a family of dimension at least $3$ of complete base-point free $g^1_6$'s whose residual is
also base-point free.}
\end{rem}

\section{Some results on Enriques surfaces}

\noindent We will use the following well-known

\begin{defn}
\label{phi}
Let $L$ be a line bundle on an Enriques surface $S$ such that $L^2 > 0$. Following \cite{cd}
we define
\[ \phi(L) =\inf \{|F.L| \; : \; F \in \Pic S, F^2 = 0, F \not\eqv 0\}. \]
\end{defn} 

\noindent This function has two important properties:
\begin{itemize} 
\item[(i)] $\phi(L)^2 \leq L^2$ (\cite[Cor.2.7.1]{cd});
\item[(ii)] If $L$ is nef, then there exists a genus one pencil $|2E|$ such that $E.L = \phi(L)$
(\cite[2.11]{co} or by \cite[Cor.2.7.1, Prop.2.7.1 and Thm.3.2.1]{cd}).
\end{itemize}

\noindent We will often use the
\begin{defn}
Let $S$ be an Enriques surface. A {\bf nodal} curve on $S$ is a smooth rational curve contained in $S$.
\end{defn} 

\noindent We will now briefly recall some results on line bundles on Enriques surfaces, proved in
\cite{klvan} and \cite{kl1}, that we will often use.

\begin{lemma} \cite[Lemma 2.2]{kl1}
\label{A}
Let $L > 0$ and  $\Delta > 0$ be divisors on an Enriques surface $S$ with $L^2 \geq 0$, $\Delta^2 =
-2$ and $k := - \Delta.L > 0$. Then there exists an $A > 0$ such that $A^2 = L^2$, $A.\Delta = k$ and 
$L \sim A + k \Delta$. Moreover if $L$ is primitive then so is $A$. 
\end{lemma}

\begin{lemma} \cite[Lemma 2.3]{kl1}
\label{lemma:nefred} 
Let $S$ be an Enriques surface and let $L$ be a line bundle on $S$ such that $L > 0, L^2 > 0$.  Let
$F > 0$ be a divisor on $S$ such that $F^2 = 0$ and $\phi(L) = |F.L|$. Then
\begin{itemize}
\item[(a)] $F.L > 0$;
\item[(b)] if $\alpha > 0$ is such that $(L - \alpha F)^2 \geq 0$, then $L - \alpha F > 0$.
\end{itemize}
\end{lemma}

\begin{lemma} \cite[Lemma 2.1]{klvan}
\label{lemma10}
Let $X$ be a smooth surface and let $A > 0$ and $B > 0$ be divisors on $X$ such that
$A^2 \geq 0$ and $B^2 \geq 0$. Then $A.B \geq 0$ with equality if and only if there exists a primitive divisor $F
> 0$ and integers $a \geq 1, b \geq 1$ such that $F^2 = 0$ and $A \eqv aF, B \eqv bF$.
\end{lemma}

\begin{defn} 
\label{def:qnef}
An effective line bundle $L$ on a K3 or Enriques surface is said to be {\bf quasi-nef} if $L^2
\geq 0$ and $L.\Delta \geq -1$ for every $\Delta$ such that $\Delta > 0$ and $\Delta^2 = -2$. 
\end{defn}

\begin{thm} \cite[Corollary 2.5]{klvan}
\label{cor:qnef} 
An effective line bundle $L$ on a K3 or Enriques surface is quasi-nef if and only if $L^2 \geq 0$
and either $h^1(L) = 0$ or $L \eqv nE$ for some $n \geq 2$ and some primitive and nef divisor $E > 0$
with $E^2 = 0$.
\end{thm}

\begin{thm} \cite[Corollary 1]{kl1}
\label{corkl1}
Let $|L|$ be a base-component free linear system on an Enriques surface $S$ such that $L^2 > 0$ and 
let $C \in |L|$ be a general curve. Then 
\[ \gon(C) = 2\phi(L) \]
unless $L$ is of one of the following types:
\begin{itemize}
\item[(a)] $L^2 = \phi(L)^2$ with $\phi(L) \geq 2$ and even. In these cases $\gon(C) = 2\phi(L) - 2$.
\item[(b)] $L^2 = \phi(L)^2 + \phi(L) - 2$ with $\phi(L) \geq 3$, $L \not\eqv 2D$ for $D$ such that 
$D^2 = 10$, $\phi(D) = 3$. In these cases $\gon(C) = 2\phi(L) - 1$ except for $\phi(L) = 3, 4$ when
$\gon(C) = 2\phi(L) - 2$.
\item[(c)] $(L^2, \phi(L)) = (30, 5)$, $(22, 4)$, $(20, 4)$, $(14, 3)$, $(12, 3)$ and $(6, 2)$. In 
these cases $\gon(C) = \lfloor \frac{L^2}{4} \rfloor + 2 = 2\phi(L) - 1$.
\end{itemize}
\end{thm}

\section{Tetragonal curves on Enriques surfaces and on surfaces \\ of degree $g - 1$ in $\PP^{g - 1}$}
\label{onenr}

\noindent Let $C$ be a smooth irreducible tetragonal curve of genus $g \geq 6$ and let $M$  be a line bundle on 
$C$ such that $H^1(M) = 0$ and $\mu_{M, \omega_C}$ is surjective. To have the surjectivity of the Gaussian map
$\Phi_{M, \omega_C}$, it is necessary, by Proposition \ref{tetragonal}(ii), that $h^0(2K_C - M - b_{2, A} A)
= 0$ for every $g^1_4$ on $C$. On the other hand when $h^0(2K_C - M) = 1$ we need that $b_{2, A} \geq 1$ for
every $g^1_4$ on $C$, that is (see \ref{b2}) $b_2(C) \geq 1$, because in this case, by Proposition
\ref{tetragonal}(ii), $h^0(2K_C - M - b_{2, A} A) = \cork \Phi_{M, \omega_C}$ is independent of $A$. As we have
seen in \ref{b2}, in the canonical embedding, $C = Y_A \cap Z_A \subset \PP^{g - 1}$ where $Y_A$ is a surface of
degree $g - 1 + b_{2, A}$ by Lemma \ref{N2perY}. Moreover $b_{2, A} = 0$ if and only if $C$ is a quadric section
of $Y_A$. Therefore saying that $b_2(C) \geq 1$ is equivalent to saying that $C$, in its canonical embedding,
can never be a quadric section of a surface $Y_A$ of degree $g - 1$ in $\PP^{g - 1}$. 

\noindent The present section we will be devoted to proving that tetragonal curves of genus $g \geq 7$, lying on
an Enriques surface and general in their linear system, in their canonical embedding, can never be a quadric
section of a surface $Y_A$ of degree $g - 1$ in $\PP^{g - 1}$. The latter fact will be then used to prove
surjectivity of Gaussian maps for such curves in our main theorem.

\noindent We start by observing that we cannot do better in genus $6$. Let $C$ be a smooth irreducible tetragonal
curve of genus $6$ and let $A$ be a $g^1_4$ on $C$. Now $K_C - A$ is a $g^2_6$ and has a base point if and only 
if $C$ is isomorphic to a plane quintic. Therefore if $C$ is not isomorphic to a plane quintic, then it has
complete base-point free $g^2_6$ and either $C$ is bielliptic or the $g^2_6$ is birational. In the latter case
the image of $C$ by the $g^2_6$ cannot have points of multiplicity higher than $2$, therefore $C$ does lie on 
$X = \Sigma_4$ and is linearly equivalent to $-2K_X$. 

\noindent Hence we can  restrict our attention to curves of genus $g \geq 7$.

\vskip .3cm
 
\noindent We will henceforth let $S$ be an Enriques surface. 

\vskip .3cm

\noindent Consider a base-point free line bundle $L$ on $S$ with $L^2 \geq 12$ and let $C \in |L|$ be a
general curve. By Theorem \ref{corkl1} we have that $C$ is not trigonal and moreover $C$ is tetragonal if
and only if $\phi(L) = 2$. 

\noindent Now assume that $\phi(L) = 2$. We have

\begin{thm} 
\label{nodelp}
Let $L$ be a base-point free line bundle on an Enriques surface with $L^2 \geq 12$ and $\phi(L) = 2$.
Then $b_2(C) \geq 1$ for a general curve $C \in |L|$. 
\end{thm}

\noindent The proof of this theorem will be essentially divided in two parts, namely a careful study of the cases
$L^2 = 12, 14$ and $16$ and an application of previous results for $L^2 \geq 18$. In both parts we will
employ the following

\begin{gerem}
\label{genrem}
{\rm Let $C$ be a tetragonal curve of genus $g$ and let $A$ be a $g^1_4$ on $C$ such that $b_{2, A} = 0$. Then, by
\ref{b2} and Lemma \ref{N2perY}, in its canonical embedding, $C$ is a quadric section of a surface $Y_A \subset
\PP^{g - 1}$ of degree $g - 1$ whence, by Remark \ref{delp2}, $C$ is contained in a surface $X$ that is either
$\Sigma_{10 - g}$, or $Bl_V C_{g - 1}$, or a smooth quadric in $\PP^3$ or $Bl_V Q$ where $Q$ is a quadric cone
in $\PP^3$, or $\PP^2$. Also $C$ is either bielliptic (in the case of $C_{g - 1}$) or linearly equivalent to $-
2K_X$.}
\end{gerem}

\noindent We start with the cases of genus $7, 8$ and $9$.

\subsection{Curves of genus $7$}

\noindent We will need the ensuing

\begin{lemma}
\label{lemmag=7}
Let $L$ be a base-point free line bundle on an Enriques surface with $L^2 = 12$ and $\phi(L) = 2$.
Let $|2E|$ be a genus one pencil such that $E.L = 2$. Then there exists a primitive divisor $E_1$
such that $E_1 > 0$, $E_1^2 = 0$, $E + E_1$ is nef, $h^0(E_1) = h^0(E_1 + K_S) = 1$ and one of the following cases occurs:
\begin{itemize}
\item[(i)] $\phi(L - 2 E) = 1$ and $L \sim 3E + 2E_1$, $E.E_1 = 1$;
\item[(ii)] $\phi(L - 2 E) = 2$ and $L \sim 3E + E_1$, $E.E_1 = 2$.
\end{itemize}
Moreover, in case (ii), for any smooth curve $C \in |L|$, we have that $h^0((E_1)_{|C}) = h^0((E_1 + K_S)_{|C}) = 2$.
\end{lemma}

\begin{proof}
We have $(L - 3E)^2 = 0$, $E.(L - 3E) = 2$ and by Lemma \ref{lemma:nefred} we can write $L \sim 3E + E'_1$ with $E'_1
> 0$, $(E'_1)^2 = 0$ and $E.E'_1 = 2$. Also $1 \leq \phi(L - 2 E) \leq \sqrt{(L - 2 E)^2} = 2$. 

\noindent If $\phi(L - 2 E) = 2$ we set $E_1 = E'_1$. Then certainly $E_1$ is primitive and we have $L \sim 3E + E_1$,
$E.E_1 = 2$, as in (ii). 

\noindent If $\phi(L - 2 E) = \phi(E + E'_1) = 1$ let $F > 0$ be a divisor such that $F^2 = 0$ and $F.(E + E'_1) = 1$ 
($F$ exists by Lemma \ref{lemma:nefred}). Then necessarily $F.E = 1, F.E'_1 = 0$ therefore $E'_1 \eqv 2F$ by Lemma
\ref{lemma10} and we can set $E_1 = F$. Replacing, if necessary, $E$ with $E + K_S$, we have that $E_1$ is primitive and
$L \sim 3E + 2E_1$, $E.E_1 = 1$, as in (i). 

\noindent Since $E_1$ is primitive, to see, in both cases (i) and (ii), that $h^0(E_1) = h^0(E_1 + K_S) = 1$, by
\cite[Cor.2.5]{klvan}, we just need to show that $E_1$ is quasi-nef. Let $\Delta > 0$ be a divisor such that $\Delta^2 =
-2$ and $k := - E_1.\Delta \geq 1$. By \cite[Lemma2.2]{kl1} we can write $E_1 \sim A + k\Delta$ for some $A > 0$ primitive
with $A^2 = 0$, $A.\Delta = k$. Now $0 \leq L.\Delta = 3E.\Delta + E'_1.\Delta \leq 3E.\Delta - 1$ gives $E.\Delta \geq
1$. From $2 \geq E.E_1 = E.A + k E.\Delta$ we get that either $k = 1$ or $k = 2, E.\Delta = 1$ and $E.A = 0$. In the
latter case we have that $E \eqv A$ by Lemma \ref{lemma10} and this is a contradiction since $A.\Delta = 2$.

\noindent Therefore we have proved that $E_1$ is quasi-nef and if $E_1.\Delta \leq -1$ then $E_1.\Delta = -1$, $E.\Delta
\geq 1$. This of course implies that $E + E_1$ is nef.

\noindent Suppose now that we are in case (ii), let $F \eqv E_1$ and let $C \in |L|$ be a smooth curve. From
the exact sequence
\[0 \hpil F - L \hpil F \hpil F_{|C} \hpil 0 \]
and the fact just proved that $h^0(F) = 1$, $h^1(F) = 0$, we see that $h^0(F_{|C}) = 1 + h^1(F - L) = 2$ since $F - L
\eqv - 3E$.
\end{proof}

\noindent The above lemma allows to exclude quickly the bielliptic case.

\begin{rem}
\label{bielliptic} 
{\rm Let $L$ be a base-point free line bundle on an Enriques surface $S$ with $L^2 = 12$ and 
$\phi(L) = 2$. Let $|2E|$ be a genus one pencil such that $E.L = 2$. Let $C$ be a general curve in $|L|$. If $b_2(C) = 0$
we can certainly say that $C$ is not bielliptic since if $A$ is a complete base-point free $g^1_4$ on $C$ we have, by
Proposition \ref{g1kondelp}(g), that $A \sim (f_1 + f_2)_{|C}$ therefore $|K_C - A| = |(f_1' + \ldots + f'_4)_{|C}|$ is
not birational. On the other hand on the Enriques surface $S$, if we pick $A = (2E)_{|C}$, using the notation of Lemma
\ref{lemmag=7}, we have that either $K_C - A \sim (E + E_1 + K_S)_{|C}$ or $K_C - A \sim (E + 2E_1 + K_S)_{|C}$. Since
the linear systems $|E + E_1 + K_S|$ and $|E + 2E_1 + K_S|$ define a map whose general fiber is finite by \cite[Thm.4.6.3
and Thm.4.5.1]{cd}, we get that $|K_C - A|$ is birational for general $C$ since $|L|$ is birational by \cite[Thm.4.6.3
and Prop.4.7.1]{cd}.}
\end{rem} 

\noindent According to the two cases in Lemma \ref{lemmag=7} we will have two propositions.

\begin{prop} 
\label{g=7}
Let $L$ be a base-point free line bundle on an Enriques surface with $L^2 = 12$ and $\phi(L) = 2$.
Let $|2E|$ be a genus one pencil such that $E.L = 2$ and suppose that $\phi(L - 2 E) = 1$.

\noindent Then $b_2(C) \geq 1$ for a general curve $C \in |L|$. 
\end{prop}
 
\begin{proof}
We use the notation of Lemma \ref{lemmag=7}. 

\noindent First we prove that either $(E + E_1)_{|C}$ or $(E + E_1 + K_S)_{|C}$ is a complete base-point
free $g^1_5$ on $C$.

\noindent To this end note that since $(E + E_1)^2 = 2$ and $E + E_1$ is nef by Lemma \ref{lemmag=7}, we
have by \cite[Prop.3.1.6 and Cor.3.1.4]{cd} that either $E + E_1$ or $E + E_1 + K_S$ is base-component free with
two base points.  Let $B \eqv E + E_1$ be the line bundle that is base-component free. As $C$ is general in
$|L|$ we have that $B_{|C}$ is base-point free. Now the exact sequence
\[ 0 \hpil B - C \hpil B \hpil B_{|C} \hpil  0 \]
shows that also $B_{|C}$ is a complete $g^1_5$ since $B - C \eqv -2 E - E_1$ whence $h^1(B - C) = 0$ because
$2 E + E_1$ is nef by Lemma \ref{lemmag=7}.

\noindent Now suppose that there exists a line bundle $A$ that is a $g^1_4$ on $C$ and is such that 
$b_{2, A} = 0$. By the general remark \ref{genrem} we know that $C$ lies on a surface $X$ (obtained by desingularizing
$Y_A$, if necessary) and either $X = \Sigma_3$, $C \sim -2K_X$ or $X = Bl_{V} C_6$ and $C$ is bielliptic. As $C$ has a
complete base-point free $g^1_5$ the second case is excluded by Proposition \ref{g1kondelp}(g) (or by Remark
\ref{bielliptic}). When $X = \Sigma_3$ by Proposition \ref{g1kondelp}(e) we know that there is a point $P \in C$ such that
either $B_{|C} \sim \widetilde{H}_{|C} - P$ or $B_{|C} \sim (2\widetilde{H} - G_1 - G_2 - G_3)_{|C} - P$. 

\noindent If $B_{|C} \sim \widetilde{H}_{|C} - P$ then 
\[ K_C \sim (3 \widetilde{H} - G_1 - G_2 - G_3)_{|C} \sim (L + K_S - B)_{|C} + \widetilde{H}_{|C} - P \]
whence
\begin{equation}
\label{eq3}
(L + K_S - B)_{|C} - P \sim (2\widetilde{H} - G_1 - G_2 - G_3)_{|C} \ \mbox{is a} \ g^2_6 \ \mbox{on C.}
\end{equation} 
If $B_{|C} \sim (2\widetilde{H} - G_1 - G_2 - G_3)_{|C} - P$ then 
\[ K_C \sim (3 \widetilde{H} - G_1 - G_2 - G_3)_{|C} \sim (L + K_S - B)_{|C} + (2\widetilde{H} - G_1 - G_2 -
G_3)_{|C} - P \] 
whence
\begin{equation}
\label{eq4}
(L + K_S - B)_{|C} - P \sim \widetilde{H}_{|C} \ \mbox{is a} \ g^2_6 \ \mbox{on C.}
\end{equation}  

\noindent But using the Enriques surface $S$ we have an exact sequence
\[ 0 \hpil K_S - B \hpil L + K_S - B \hpil (L + K_S - B)_{|C} \hpil  0 \]
and $L + K_S - B \eqv 2 E + E_1$, $h^1(K_S - B) = 0$ by Lemma \ref{lemmag=7}, whence $(L + K_S - B)_{|C}$ is 
a base-point free $g^2_7$ on $C$, contradicting (\ref{eq3}) and (\ref{eq4}).
\end{proof}

\noindent Now the other case. 

\begin{prop}
\label{g^2_6particolari-tris}
Let $L$ be a base-point free line bundle on an Enriques surface $S$ with $L^2 = 12$ and $\phi(L) = 2$.
Let $|2E|$ be a genus one pencil such that $E.L = 2$ and suppose that $\phi(L - 2 E) = 2$.
Then the general curve in $|L|$ possesses no $g^2_6$ and satisfies $b_2(C) \geq 1$.
\end{prop}

\begin{proof}
The proof will be a variant of the method of \cite[Section4]{kl1}. By Lemma \ref{lemmag=7} we have $L \sim 3E
+ E_1$ with $E > 0$, $E_1 > 0$ both primitive, $E^2 = E_1^2 = 0$, $E$ and $E + E_1$ are nef and $E.E_1 = 2$. Let $D =
2E + E_1$ so that $D^2 = 8$, $\phi(D) = 2$, $D.L = 10$ and $D$ is nef, whence base-point free by
\cite[Prop.3.1.6, Prop.3.1.4 and Thm.4.4.1 ]{cd}. 

\noindent Now recall that by \cite[Thm.4.6.3 and Thm.4.7.1]{cd} the linear system $|D|$ defines a birational morphism
$\varphi_D : S \to \overline{S} \subset \PP^4$ onto a surface $\overline{S}$ having some rational double points,
corresponding to nodal curves $R \subset S$ such that $D.R = 0$, and two double lines, namely $\varphi_D(E)$ and
$\varphi_D(E + K_S)$. More precisely by \cite[Prop.3.7]{kn2} we see that if $Z \subset S$ is any zero-dimensional
subscheme of length two not imposing independent conditions to $|D|$ then either $Z \subset E$ or $Z \subset E + K_S$ or
any point $x \in \Supp(Z)$ lies on some nodal curve contracted by $\varphi_D$. Observe that if $R \subset S$ is a nodal
curve contracted by $\varphi_D$, then $0 = D.R = E.R + (E + E_1).R$ whence $E.R = E_1.R = 0$ by the nefness of $E$ and of
$E + E_1$. This implies that $C.R = 0$, whence that $C \cap R = \emptyset$, for any $C \in |L|_{sm}$. Also, if
$\overline{S}$ contains a line different from the two double lines, then this line is image of a nodal curve $\Gamma
\subset S$ such that $D.\Gamma = 1$ whence, using again the nefness of $E + E_1$, we have that either $E.\Gamma = 0$,
$E_1.\Gamma = 1$ or $E.\Gamma = 1$, $E_1.\Gamma = -1$. This implies that $C.\Gamma = 1, 2$ for any $C \in |L|_{sm}$. In
particular, since $C.E = 2$, we find that for each line on $\overline{S}$ its inverse image in $S$ can contain at most two
points of any $C \in |L|_{sm}$. Moreover $\overline{S}$ contains finitely many lines, namely the two lines $\varphi_D(E)$,
$\varphi_D(E + K_S)$ and the images of the finitely many irreducible curves $\Gamma \subset S$ such that $D.\Gamma = 1$
(these are finitely many since if $D.\Gamma = 1$ we get $\Gamma^2 = -2$).

\noindent By Remark \ref{bielliptic} we know that there is a proper closed subset $B \subset |L|_{sm}$ such that every
element in $B$ is bielliptic and by Theorem \ref{corkl1} there is another proper closed subset $B_3 \subset |L|_{sm}$
such that every element in $B_3$ is trigonal or hyperelliptic and any element of $\U := |L|_{sm} - (B \cup B_3)$ is
tetragonal. We set $B^2_6$ for the closed subset of $|L|_{sm}$ whose elements correspond to curves having a $g^2_6$. 

\noindent The goal will be to prove that the open subset $|L|_{sm} - (B \cup B_3 \cup B^2_6)$ is nonempty. 

\noindent We will therefore suppose that it is empty, so that every $C \in \U$ has a linear series $A_C$ that is a
$g^2_6$ on $C$.

\noindent Since $h^0(D_{|C} - A_C) = h^0(\omega_C - A_C - (E + K_S)_{|C}) \geq h^1(A_C) - 2 \geq 1$, we see that
there exists an effective divisor $T$ of degree $4$ on $C$ such that $T \sim D_{|C} - A_C$. 

\begin{claim}
\label{span}
For each $T$ as above we have $h^0(\I_{T/S}(D)) = 3$ and $h^0(\I_{T/S}(L)) = 4$.
\end{claim}
\begin{proof}
The first part of the claim follows by the exact sequence
\begin{equation}
\label{sp}
0 \hpil \O_S(- E) \hpil \I_{T/S}(D) \hpil \I_{T/C}(D) \hpil 0
\end{equation}
since then $h^0(\I_{T/S}(D)) = h^0(\I_{T/C}(D)) = h^0(A_C) = 3$.

\noindent To see the second part of the claim consider the exact sequence
\[ 0 \hpil \O_S \hpil \I_{T/S} \otimes L \hpil \I_{T/C} \otimes L \hpil 0 \]
so that $h^0(\I_{T/S} \otimes L) = 1 + h^0(L_{|C} - T) = 1 + h^0(A_C + E_{|C})$. 

\noindent We will prove that $h^0(A_C + E_{|C}) = 3$. Now $h^0(A_C + E_{|C}) \geq h^0(A_C) = 3$ and we need to exclude
that $h^0(A_C + E_{|C}) \geq 4$. 

\noindent Assume henceforth that $h^0(A_C + E_{|C}) \geq 4$. 

\noindent Since $\deg(A_C + E_{|C}) = 8$ and $\Cliff(C) = 2$, if $h^0(A_C + E_{|C}) \geq 4$, we must have $h^0(A_C +
E_{|C}) = 4$, therefore $h^0(\I_{T/S} \otimes L) = 5$. Since $h^0(L) = 7$ we see that there is a zero-dimensional
subscheme $Z \subset T$ such that $\length(Z) = 3$ and $h^0(\I_{Z/S} \otimes L) = 5$. We claim that there is a proper
subscheme $Z' \subset Z$ such that $\length(Z') = 2$ and $h^0(\I_{Z'/S} \otimes L) \geq 6$. In fact if for every proper
subscheme $Z' \subset Z$ with $\length(Z') = 2$ we have $h^0(\I_{Z'/S} \otimes L) = 5$ then $Z$ is in special position
with respect to $L + K_S$ and, since $L^2 = 4 \length(Z) = 12$, we deduce by \cite[Prop.3.7]{kn2} that there is an
effective divisor $B$ such that $Z \subset B$ and $L.B \leq B^2 + 3 \leq 6$. Since $B.L \geq 3$ we get that
\begin{equation}
\label{bog}
3 \leq L.B \leq B^2 + 3 \leq 6 
\end{equation}
whence $0 \leq B^2 \leq 2$. Note that for any $F > 0$ with $F^2 = 0$ we have either $F.L \geq 4$ or $F \eqv E$
(whence $F.L = 2$). Now if $B^2 = 2$ we can write $B \sim F_1 + F_2$ with $F_i > 0$, $F_i^2 = 0$ for $i = 1, 2$ and
$F_1.F_2 = 1$. By (\ref{bog}) we have $L.F_1 + L.F_2 = L.B \leq 5$, whence the contradiction $F_1 \eqv E \eqv F_2$.
Therefore $B^2 = 0$ and $L.B = 3$ by (\ref{bog}), again a contradiction.

\noindent We have therefore proved that there is a proper subscheme $Z' \subset Z \subset C$ such that $\length(Z') = 2$
and $h^0(\I_{Z'/S} \otimes L) \geq 6$, whence $h^0(\I_{Z'/S} \otimes L) = 6$ as $L$ is base-point free and therefore
$Z'$ is not separated by the morphism $\varphi_L : S \to \PP^6$. Now recall that by \cite[Thm.4.6.3, Prop.4.7.1 and
Cor.1, p.283]{cd} $\varphi_L$ is a birational morphism onto a surface having some rational double points, corresponding to
nodal curves $R \subset S$ such that $L.R = 0$, and two double lines, namely $\varphi_L(E)$ and $\varphi_L(E + K_S)$ and
that $\varphi_L$ is an isomorphism outside $E, E + K_S$ and the nodal curves contracted. In particular we deduce that
either $Z' = C \cap E$ or $Z' = C \cap E'$. We claim that this implies that either $T \sim (2E)_{|C}$ or $T \sim (2E +
K_S)_{|C}$. 

\noindent To see the latter suppose for example that $Z' = C \cap E$ and set $W = T - Z'$ on $C$. Then
$\length(W) = 2$ and $4 = h^0(L_{|C} - T) = h^0(D_{|C} + Z' - T) = h^0(D_{|C} - W)$, whence the exact sequence
\[ 0 \hpil \O_S(- E) \hpil \I_{W/S} \otimes D \hpil \I_{W/C} \otimes D \hpil 0 \]
shows that $h^0(\I_{W/S} \otimes D) = 4$. Therefore $W$ is not separated by the morphism $\varphi_D : S \to \PP^4$. As 
$C \cap R = \emptyset$, for any nodal curve $R$ contracted by $\varphi_D$ we have that either $W \sim E_{|C}$ or
$W \sim (E + K_S)_{|C}$, whence either $T \sim (2E)_{|C}$ or $T \sim (2E + K_S)_{|C}$. 

\noindent Finally since we know that $T \sim D_{|C} - A_C$ we deduce that either $A_C \sim (E_1)_{|C}$ or $A_C \sim
(E_1 + K_S)_{|C}$, but this contradicts Lemma \ref{lemmag=7}.
\end{proof}

\noindent {\it Continuation of the proof of Proposition {\rm \ref{g^2_6particolari-tris}}}. 
Consider the following incidence subscheme of $\Hilb^4 (S) \times \U$ :
\[ \J = \{ (T, C) \; : \; T \in \Hilb^4 (S), C \in \U, \; T \subset C \ \mbox{and} \ h^0(D_{|C} - T) \geq 3 \} \]
together with its two projections $\pi : \J \to \Hilb^4 (S)$ and $p : \J \to \U.$

\noindent Our assumption that any $C \in \U$ carries a $g^2_6$ implies, as we have seen, that $p$ is surjective, whence we
deduce that $\J$ has an irreducible component $\J_0$ such that $\dim \J_0 \geq 6$. Since the fibers of
$\pi$ have dimension at most $h^0(\I_{T/S}(L)) - 1 = 3$ by Claim \ref{span}, we get that $\dim \pi(\J_0) \geq 3$. 

\noindent Using $\pi(\J_0)$ we build up an incidence subscheme of $\pi(\J_0) \times |D|$ :
\[ \I = \{ (T, D') \; : \; T \in \pi(\J_0), D' \in |D|, \; T \subset D' \} \]
together with its two projections 
\begin{equation}
\label{jdc}
f : \I \to |D| \ \mbox{and} \ h : \I \to \pi(\J_0). 
\end{equation}
By (\ref{sp}) and the definition of $\pi(\J_0)$ we have that $h$ is surjective. Since the fibers of $h$ have
dimension $h^0(\I_{T/S}(D)) - 1 = 2$ by Claim \ref{span}, we find that $\I$ has an irreducible component $\I_0$ such
that $\dim \I_0 \geq 5$. 

\noindent To show that this fact leads to a contradiction let us return to the morphism $\varphi_D : S \to \overline{S}
\subset \PP^4$.

\noindent A general hyperplane section $\overline{D} = \overline{S} \cap H \subset \PP^3$ is a curve of degree
$8$ with two nodes, whence of arithmetic genus $7$. Consider, for $i = 2, 3$, the exact sequence
\[0 \hpil \O_{\overline{S}}(i - 1) \hpil \O_{\overline{S}}(i) \hpil \O_{\overline{D}}(i) \hpil 0. \]
Using Riemann-Roch on $\overline{D}$ we get 
\[ h^0(\O_{\overline{S}}(3)) \leq h^0(\O_{\overline{S}}(2)) + h^0(\O_{\overline{D}}(3)) \leq h^0(\O_{\overline{S}}(1)) +
h^0(\O_{\overline{D}}(2)) + h^0(\O_{\overline{D}}(3)) = 33 \]
whence $h^0(\I_{\overline{S}/\PP^4}(3)) \geq 2$ and therefore there is a plane $\overline{P} \subset \PP^4$ such that
$\overline{S} \cup \overline{P}$ is a complete intersection of two cubics in $\PP^4$. 

\noindent Now every $T \in \pi(\J_0)$ has three important properties. First of all we know that $T \subset C$ for
some $C \in \U$ and $C \cap R = \emptyset$ for every nodal curve $R$ contracted by $\varphi_D$, therefore also
$T \cap R = \emptyset$ for every nodal curve $R$ contracted by $\varphi_D$. Secondly, since $C.E = 2$, we get
that $\length(T \cap E) \leq 2$ and $\length(T \cap (E + K_S)) \leq 2$. Thirdly the linear span $l_T := <\varphi_D(T)>
\subset \PP^4$ is a line by Claim \ref{span}. Moreover let us prove that we cannot have infinitely many elements $T \in
\pi(\J_0)$ such that $l_T$ is the same line. Suppose to the contrary that there is an infinite set $Z \subset
\pi(\J_0)$ and a line $l \subset \PP^4$ such that $l_T = l$ for every $T \in Z$. If $l$ is not contained in
$\overline{S}$ then it meets $\overline{S}$ in finitely many points, therefore there is a point $P \in l$ and an
infinite set $V \subset S$ such that $\varphi_D(x) = P$ for every $x \in V$ and each $x \in V$ lies on some $T \in Z$.
Now $V \subset \varphi_D^{-1}(P)$ therefore $\varphi_D^{-1}(P)$, being infinite, must be a nodal curve contracted by
$\varphi_D$ (recall that $\varphi_D$ is $2$ to $1$ on $E$ and $E + K_S$) and this is absurd since for any for $x \in V$
we have that $x \in T$ for some $T \in Z$ and we know that $T \cap R = \emptyset$ for every nodal curve $R$ contracted
by $\varphi_D$. Therefore $l$ is contained in $\overline{S}$ and all $T \in Z$ lie in $\varphi_D^{-1}(l) \subset S$ and
this is absurd since each $T$ is contained in some $C \in \U$ and we know that $\varphi_D^{-1}(l)$ can contain at most
two points of any $C \in \U$.

\noindent Since $\dim \pi(\J_0) \geq 3$ we have that there is a family of lines $l_T := <\varphi_D(T)>$ of
dimension at least $3$ meeting $\overline{S}$ along $\varphi_D(T)$.

\noindent Now let $T \in \pi(\J_0)$ be a general element. We cannot have that $\length(\varphi_D(T)) \geq 4$, else
$\varphi_D(T)$ is contained in $l_T \cap F_3$ for every cubic $F_3$ containing $\overline{S}$, that is $l_T$ is contained
in $\overline{S} \cup \overline{P}$, a contradiction since $\overline{S}$ contains finitely many lines and of course
$\overline{P} \cong \PP^2$ contains a $2$-dimensional family of lines.

\noindent Therefore $\length(\varphi_D(T)) \leq 3$ for a general $T \in \pi(\J_0)$, whence such a $T$ is not mapped
isomorphically by $\varphi_D$ and therefore it does not lie in the open subset $S - E \cup (E + K_S) \cup R_1 \cup \ldots
\cup R_n$, where $R_1, \ldots, R_n$ are the nodal curves contracted by $\varphi_D$. This means that for a general
$T \in \pi(\J_0)$ we have that $\varphi_D(T) \cap \Sing(\overline{S}) \not= \emptyset$ and therefore also $\varphi_D(T)
\cap \Sing(\overline{D_0}) \not= \emptyset$ for any $D_0$ containing $T$.

\noindent Now consider the map $f_0 := f_{|\I_0} : \I_0 \to |D|$ from (\ref{jdc}). Certainly $f_0$ cannot be surjective,
for otherwise a general hyperplane section $\overline{D_0}$ of $\overline{S}$ would have infinitely many lines $l_T$
passing through a fixed node of $\overline{D_0}$. Hence the projection of $\overline{D_0}$ from that node would give
either a $2$ to $1$ map of $D_0$ onto a singular (since $\overline{D_0}$ has two nodes) plane cubic, whence $D_0$ would be
hyperelliptic, or a $3$ to $1$ map onto a conic, whence $D_0$ would be trigonal. Therefore $\gon(D_0) \leq 3$ for a
general $D_0 \in |D|$, contradicting Theorem \ref{corkl1}.

\noindent Hence $\dim f_0(\I_0) \leq 3$ and let $D_0 \in f_0(\I_0)$ be a general element. If $\dim f_0(\I_0)
\leq 2$ we get that $\dim f_0^{-1}(D_0) \geq 3$, whence a general element $(T, D_0) \in f_0^{-1}(D_0)$ is such that at
least three points of $T$ are general on $D_0$. But this is a contradiction since these points give rise to three general
points of $\overline{D_0}$ that span a line.

\noindent Therefore $\dim f_0(\I_0) = 3$. We will first prove that this implies that $D_0$ is reducible. 

\noindent Suppose that $D_0$ is irreducible. Note that both $\overline{D_0}$ and $D_0$ are reduced, for otherwise we would
have that $D_0 \sim m\Delta$ for some $\Delta > 0$ and some $m \geq 2$, but then $D_0^2 = 8$ implies $m = 2$, whence that
$E_1$ is $2$-divisible, a contradiction. Since $\dim \I_0 \geq 5$ we have that $\dim f_0^{-1}(D_0) \geq 2$ and for each
$(T, D_0) \in f_0^{-1}(D_0)$ we know that $l_T = <\varphi_D(T)> \subset \PP^3$ is a line. Moreover we showed above that
we cannot have infinitely many divisors $T$'s such that $l_T \subset \PP^3$ is the same line, therefore $\overline{D_0}$
has a family of dimension at least two of lines $l_T$ meeting $\overline{D_0}$ along $\varphi_D(T)$ and
$\length(\varphi_D(T)) \leq 3$ for a general such $T$. Hence we have a family of dimension at least two of lines meeting
$\overline{D_0}$ on a singular point $P_0$ of $\overline{D_0}$ and meeting it furthermore at two points (possibly
coinciding). Also, since $T \cap R = \emptyset$ for every nodal curve $R \subset S$ contracted by $\varphi_D$, we see
that $\varphi_D(T)$ is not a point. Now a general projection $D_0'$ of $\overline{D_0}$ in $\PP^2$ has the same property,
namely that the general secant line to $D_0'$ goes through a fixed point (the projection of $P_0$) and this is absurd
since $D_0'$ is not a line. This proves that $D_0$ is reducible and we can now assume that
\[f_0(\I_0) \subset \{D' \in |D| : D' \ \rm{is \ reducible} \}. \]

\noindent To exclude this case we will therefore study the reducible locus of $|D|$. To this end we first prove
the following two facts.

\begin{claim}
\label{nomoving}
There is no decomposition $D \sim A + B$ with $h^0(A) \geq 2$ and $h^0(B) \geq 2$.
\end{claim}
\begin{proof}
Suppose such a decomposition exists. Then we get $A.D \geq 2\phi(D) = 4$ and similarly $B.D \geq 4$, whence $A.D = B.D =
4$, since $D^2 = 8$. Let $A \sim F_A + M_A$, $B \sim F_B + M_B$ be the decompositions into base-components and moving
parts of $|A|$ and $|B|$. Then $h^0(M_A) \geq 2$ and $h^0(M_B) \geq 2$, whence, as above, $M_A.D = M_B.D = 4$. Now by
\cite[Prop.3.1.4]{cd} either $M_A \sim 2hE'$ for some genus one pencil $|2E'|$ or $M_A^2 > 0$. In both cases we can
write $M_A \sim \sum\limits_{i = 1}^n F_i$ with $F_i > 0$, $F_i^2 = 0$ and $n \geq 2$, therefore $4 = M_A.D \geq
n\phi(D) = 2n$. Hence $n = 2$ and $M_A \sim 2E$, since for any $F > 0$ with $F^2 = 0$ and $F.D = 2$ we must have
$F \eqv E$. Similarly $M_B \sim 2E$ and therefore $2E + E_1 = D \geq 4 E$. But then $E_1 \geq 2E$ whence $h^0(E_1)
\geq 2$, a contradiction by Lemma \ref{lemmag=7}.
\end{proof}

\begin{claim}
\label{contracted}
Let $D \sim \Delta + M$ for some $\Delta > 0$ and $M > 0$ with $M^2 \geq 6$. Then $M^2 = 6$, $\Delta^2 = -2$,
$D.\Delta = 0$.
\end{claim}

\begin{proof}
By Riemann-Roch we have that $h^0(M) \geq 4$, whence, by Claim \ref{nomoving}, $h^0(\Delta) = 1$. Hence
$\Delta^2 \leq 0$ by Riemann-Roch and $M^2 = (D - \Delta)^2 = 8 + \Delta^2 - 2 D.\Delta \geq 6$, so that 
\[ 2 D.\Delta \leq 2 + \Delta^2. \]
If $\Delta^2 = 0$ we find the contradiction $2 \geq 2 D.\Delta \geq 2 \phi(D) = 4$. If $\Delta^2 \leq -2$, by the
nefness of $D$, we find that $0 \geq 2 D.\Delta \geq 0$, that is $M^2 = 6$, $\Delta^2 = -2$ and $D.\Delta = 0$.
\end{proof}

\noindent Now the reducible locus:

\begin{claim}
\label{redlocus}
Let $W$ be an irreducible subvariety of $\{D' \in |D| : D' \ \rm{is \ reducible} \}$ such that $\dim W = 3$. Then 
there is a divisor $G_W > 0$ with $h^0(G_W) = 1$ and such that if $M \sim D - G_W$ then $|M|$ is base-component free and
every curve $D' \in W$ is $D' = G_W + M'$ for some $M' \in |M|$. Moreover $M^2 = 6$, $G_W^2 = -2$ and $D.G_W = 0$.
\end{claim}

\begin{proof}
Let $D'$ be an element of $W$. Since $D'$ is reducible we have that $D' = G + B$ with $G > 0$, $B > 0$ and, by
Claim \ref{nomoving}, we can assume that $h^0(G) = 1$. Since the divisor classes $G > 0$ such that $D - G > 0$ are
finitely many, we see that $h^0(B) \geq 4$. Let $G'$ be the base component of $|B|$ and let $M$ be its moving part.
Then also $h^0(M) \geq 4$ and $M^2 > 0$, for otherwise we have $M^2 = 0$ whence, by \cite[Prop.3.1.4]{cd}, we get that
$M \sim 2hE'$, with $|2E'|$ a genus one pencil and $h + 1 = h^0(M) \geq 4$, contradicting Claim \ref{nomoving}, since
then $D \sim 2E' + 2(h - 1) E' + G + G'$ and $h^0(2E') = 2$, $h^0(2(h - 1) E' + G + G') \geq 2$. Therefore $1 +
\frac{M^2}{2} = h^0(M) \geq 4$, whence $M^2 \geq 6$ and of course $D \sim G + G' + M$ with $G + G' > 0$ and $h^0(G +
G') = 1$ by Claim \ref{nomoving}. By Claim \ref{contracted} we have that $M^2 = 6$, $(G + G')^2 =
-2$ and $D.(G + G') = 0$. Therefore $h^0(B) = h^0(M) = 4$.

\noindent Since the possible $G + G'$ are finitely many, we get that $\dim W = \dim |M| = 3$. Let $G_1, \ldots, G_n$
be the finite set of  divisors $G > 0$ such that $D - G > 0$ and let $B_i = D - G_i$ for $i = 1, \ldots, n$. We have
seen that for every $D' \in W$ there is an $i \in \{1, \ldots, n\}$ and a divisor $B' \in |B_i|$ so that $D' = G_i +
B'$. Let $\phi_i : |B_i| \to |D|$ be the natural inclusion defined by $\phi_i(B) = B + G_i$. Then
\[ W \subset \bigcup_{i=1}^n \Im \phi_i \]
and since $\Im \phi_i \cong |B_i|$ is a closed subset of $|D|$ and $W$ is irreducible, we deduce that there is some
$G_W$ with $h^0(G_W) = 1$, $D' \sim G_W + M$ and every curve $D' \in W$ is $D' = G_W + M'$ for some $M' \in |M|$. 
Finally the remaining part follows by Claim \ref{contracted}.
\end{proof}

\noindent {\it Conclusion of the proof of Proposition {\rm \ref{g^2_6particolari-tris}}}.
Recall that $\dim f_0(\I_0) = 3$ and that a general element $D_0 \in f_0(\I_0)$ is reducible. By Claim \ref{redlocus},
there is a $G > 0$ with $h^0(G) = 1$ and such that if $M \sim D - G$ then $|M|$ is base-component free, $M^2 = 6$, $G^2 =
-2$, $D.G = 0$ and every curve $D' \in f_0(\I_0)$ is $D' = G + M'$ for some  $M' \in |M|$. Moreover note that every
irreducible component of $G$ is a nodal curve contracted by $\varphi_D$.

\noindent Therefore $D_0 = \bigcup\limits_{i = 1}^n R_i \cup M_0$ where the $R_i$'s are nodal curves contracted by
$\varphi_D$ and $M_0$ is general in $|M|$. Now $M_0$ is a smooth irreducible curve by \cite[Prop.3.1.4 and
Thm.4.10.2]{cd} and $\varphi_D(M_0)$ is a nondegenerate (since $h^0(\sum\limits_{i = 1}^n R_i) = 1$) integral curve in
$\PP^3$. On the other hand we know that on $D_0$ there is a family of dimension at least $2$ of divisors $T$ such that
$(T, D_0) \in f_0^{-1}(D_0)$ and each $T$ gets mapped to a line $l_T$ by $\varphi_D$. Since for each $T$ we have that $T
\cap R_i = 0$ for all $i = 1, \ldots, n$, we deduce that all these $T$'s lie in $M_0$ and this gives a contradiction
since then $\varphi_D(M_0)$ would have a two dimensional family of lines $l_T$ as above. 

\noindent We have therefore proved that the general curve $C \in |L|$ possesses no $g^2_6$. 

\noindent To see that it satisfies $b_2(C) \geq 1$ suppose that there exists a line bundle $A$ that is a $g^1_4$ on
$C$ and is such that $b_{2, A} = 0$. By the general remark \ref{genrem} we know that $C$ lies on a surface $X$ (obtained
by desingularizing $Y_A$, if necessary) and either $X = \Sigma_3$, $C \sim -2K_X$ or $X = Bl_{V} C_6$ and $C$ is
bielliptic. But this is clearly a contradiction since in both cases $C$ carries $g^2_6$'s.
\end{proof}
\renewcommand{\proofname}{Proof}

\subsection{Curves of genus $8$}

\begin{prop} 
\label{g=8}
Let $L$ be a base-point free line bundle on an Enriques surface with $L^2 = 14$ and $\phi(L) = 2$.

\noindent Then $b_2(C) \geq 1$ for a general curve $C \in |L|$. 
\end{prop}

\noindent We will use the following

\begin{lemma} 
\label{lemmag=8}
Let $L$ be a base-point free line bundle on an Enriques surface with $L^2 = 14$ and $\phi(L) = 2$.
Let $|2E|$ be a genus one pencil such that $E.L = 2$. Then there exists two primitive divisors $E_1,
E_2$ such that $E_i > 0$, $E_i^2 = 0$, $E.E_i = E_1.E_2 = 1$ for $i = 1, 2$,
\[ L \sim 3E + E_1 + E_2 \]
and
\begin{itemize}
\item[(i)] $E + E_1$ is nef;
\item[(ii)] either $2E + E_2$ is nef or there exists a nodal curve $\Gamma$ such that $E_2 \eqv E_1 +
\Gamma$, $E.\Gamma = 0, E_1.\Gamma = 1, E_2.\Gamma = -1$. In particular $2E + E_2$ is quasi-nef.
\end{itemize}
Moreover let $C \in |L|$ be a general curve. Then
\begin{itemize}
\item[(iii)] either $(E + E_1)_{|C}$ or $(E + E_1 + K_S)_{|C}$ is a complete base-point 
free $g^1_6$ on $C$;
\item[(iv)] $(2E + E_2)_{|C}$ and $(2E + E_2 + K_S)_{|C}$ are complete base-point free $g^2_8$'s on $C$.
\end{itemize}
\end{lemma}

\begin{proof}
Using Lemma \ref{lemma:nefred} and Lemma \ref{lemma10} we can write $L \sim 3E + E_1 + E_2$ with
$E_i > 0$ primitive, $E_i^2 = 0$ and $E.E_i = E_1.E_2 = 1$, $i = 1, 2$.

\noindent We now claim that we can assume that $E + E_1$ is nef. 

\noindent Suppose that there is a nodal curve $\Gamma$ such that $\Gamma.(E + E_1) < 0$. Then 
$E_1.\Gamma \leq - 1 - E.\Gamma \leq -1$ and $k := - E_1.\Gamma \geq 1 + E.\Gamma \geq 1$. By Lemma \ref{A}, we can
write $E_1 \sim A + k \Gamma$ with $A > 0$ primitive with $A^2 = 0$. If $E.\Gamma > 0$ we have that
$k \geq 2$ giving the contradiction $1 = E.E_1 = E.A + k E.\Gamma \geq 2$. Therefore $E.\Gamma = 0$ and the
nefness of $L$ implies that $E_2.\Gamma > 0$. From $1 = E_2.E_1 = E_2.A + k E_2.\Gamma \geq 1$ we deduce that $k
= 1$ and $E_2.A = 0$ whence $E_2 \eqv A$ by Lemma \ref{lemma10} and therefore $E_1 \eqv E_2 + \Gamma$. Now
if in addition we have that also $E + E_2$ is not nef then the same argument above shows that there is a nodal
curve $\Gamma'$ such that $E_2 \eqv E_1 + \Gamma'$, giving the contradiction $\Gamma + \Gamma' \eqv 0$. Therefore
either $E + E_1$ or $E + E_2$ is nef and (i) is proved.

\noindent Now let $\Delta > 0$ be such that $\Delta^2 = -2$, $\Delta.(2E + E_2) < 0$. Then $E_2.\Delta \leq
- 1 - 2E.\Delta \leq -1$ and $k := - E_2.\Delta \geq 1 + 2E.\Delta \geq 1$. By Lemma \ref{A}, we can write 
$E_2 \sim A + k \Delta$ with $A > 0$ primitive with $A^2 = 0$. If $E.\Delta > 0$ we have that $k
\geq 3$ giving the contradiction $1 = E.E_2 = E.A + k E.\Delta \geq 3$. Therefore $E.\Delta = 0$ and the nefness
of $L$ implies that $E_1.\Delta > 0$. From $1 = E_1.E_2 = E_1.A + k E_1.\Delta \geq 1$ we deduce that $k = 1$
and $E_1.A = 0$, whence $E_1 \eqv A$ by Lemma \ref{lemma10} and therefore $E_2 \eqv E_1 + \Delta$. Hence $2E
+ E_2$ is quasi-nef and if it is not nef then we can choose $\Delta$ to be a nodal curve. This proves (ii).

\noindent To see (iii) note that since $(E + E_1)^2 = 2$ and $E + E_1$ is nef by (i), we
have by \cite[Prop.3.1.6 and Cor.3.1.4]{cd} that either $E + E_1$ or $E + E_1 + K_S$ is base-component free with
two base points.  Let $B \eqv E + E_1$ be the line bundle that is base-component free. As $C$ is general in
$|L|$ we have that $B_{|C}$ is base-point free. Now the exact sequence
\[ 0 \hpil B - C \hpil B \hpil B_{|C} \hpil  0 \]
shows that also $B_{|C}$ is a complete $g^1_6$ since $B - C \eqv -2 E - E_2$ whence $h^1(B - C) = 0$ by
Theorem \ref{cor:qnef} because $2 E + E_2$ is quasi-nef.

\noindent To see (iv) note that if $2E + E_2$ is nef then it is base-component free with two base points by 
\cite[Prop.3.1.6, Prop.3.1.4 and Thm.4.4.1]{cd} whence $(2E + E_2)_{|C}$ is base-point free, as $C$ is general.
The same argument shows that $(2E + E_1)_{|C}$ and $(2E + E_1 + K_S)_{|C}$ are base-point free by (i). Now if $2E
+ E_2$ is not nef then $2E + E_2 \eqv 2E + E_1 + \Gamma$ by (ii) whence again $(2E + E_2)_{|C}$ is base-point free,
since $\Gamma.C = 0$. Now the exact sequence
\[ 0 \hpil - E - E_1 \hpil 2E + E_2 \hpil (2E + E_2)_{|C} \hpil  0 \]
shows that also $(2E + E_2)_{|C}$ is a complete $g^2_8$ since $h^1(- E - E_1) = 0$ because $E + E_1$ is nef by
(i). Similarly we can show the same for $(2E + E_2 + K_S)_{|C}$.
\end{proof}

\noindent Before proving Proposition \ref{g=8} we use the above lemma to deal with the case of $\Sigma_2$. This
is used also in the proof of Proposition 4.17 in \cite{kl1}.

\begin{lemma} 
\label{lemmag=8sigma}
Let $L$ be a base-point free line bundle on an Enriques surface with $L^2 = 14$ and $\phi(L) = 2$.
Then the general curve $C \in |L|$ cannot be isomorphic to a curve linearly equivalent to $-2 K_X$ on
$X = \Sigma_2$. 
\end{lemma}

\begin{proof}
By Lemma \ref{lemmag=8}(iii) there is a line bundle $B$ such that $B \eqv E + E_1$ and $B_{|C}$ is a base-point
free complete $g^1_6$ on $C$. By Proposition \ref{g1kondelp}(c) there are two points $P_1, P_2 \in C$ such that
$B_{|C} \sim (2\widetilde{H} - G_1 - G_2)_{|C} - P_1 - P_2$. 

\noindent Now
\[ K_C \sim (3 \widetilde{H} - G_1 - G_2)_{|C} \sim B_{|C} + P_1 + P_2 + \widetilde{H}_{|C} \]  
whence
\begin{equation}
\label{g26}
(L + K_S - B)_{|C} - P_1 - P_2 \sim \widetilde{H}_{|C} \ \mbox{is a} \ g^2_6 \ \mbox{on C.}
\end{equation}

\noindent On the other hand by Lemma \ref{lemmag=8}(iv) we have that $(L + K_S - B)_{|C}$ is a base-point free
$g^2_8$ on $C$ and this contradicts (\ref{g26}).
\end{proof}

\renewcommand{\proofname}{Proof of Proposition {\rm \ref{g=8}}}  
\begin{proof}
Suppose that there exists a line bundle $A$ that is a $g^1_4$ on $C$ and is such that $b_{2, A} =
0$. By the general remark \ref{genrem} we know that $C$ lies on a surface $X$ (obtained by desingularizing $Y_A$, if
necessary) and either $X = \Sigma_2$, $C \sim -2K_X$ or $X = Bl_{V} C_7$ and $C$ is bielliptic. The latter case is
excluded since, by \cite[Prop.4.17]{kl1}, $C$ has a unique $g^1_4$ while the first case was excluded in Lemma
\ref{lemmag=8sigma}.
\end{proof}
\renewcommand{\proofname}{Proof}

\subsection{Curves of genus $9$}

\begin{prop} 
\label{g=9}
Let $L$ be a base-point free line bundle on an Enriques surface with $L^2 = 16$ and $\phi(L) = 2$.

\noindent Then $b_2(C) \geq 1$ for a general curve $C \in |L|$. 
\end{prop}

\noindent We will use the following

\begin{lemma} 
\label{lemmag=9}
Let $L$ be a base-point free line bundle on an Enriques surface with $L^2 = 16$ and $\phi(L) = 2$.
Let $|2E|$ be a genus one pencil such that $E.L = 2$. Then there exists a divisor $E_1$ such that $E_1 > 0$,
$E_1^2 = 0$, $E.E_1 = 2$ and
\[ L \sim 4E + E_1. \]
Moreover if $H^1(E_1 + K_S) \neq 0$ there exists a divisor $E_2$ such that $E_2 > 0$,
$E_2^2 = 0$, $E.E_2 = 1$, $E_1 \eqv 2E_2$ and $E + E_2$ is base-component free.
\end{lemma}

\begin{proof}
Since $(L - 4E)^2 = 0$ and $E.(L - 4 E) = 2$, by Lemma \ref{lemma:nefred} we can write $L \sim 4E + E_1$ with $E_1 > 0$,
$E_1^2 = 0$ and $E.E_1 = 2$.

\noindent By Theorem \ref{cor:qnef} if $H^1(E_1 + K_S) \neq 0$ then either $E_1 \eqv nE'$ for $n \geq 2$ and
some genus one pencil $|2E'|$ or $E_1$ is not quasi-nef. In the first case we have $2 = n E.E'$ whence $n = 2,
E.E' = 1$ and we set $E'_2 = E'$. Also $E + E'_2$ is nef in this case.

\noindent If $E_1$ is not quasi-nef there exists a $\Delta > 0$ such that $\Delta^2 = -2$, $\Delta.E_1
\leq -2$. By Lemma \ref{A}, we can write $E_1 \sim A + k \Delta$ with $A > 0$, $A^2 = 0$, $A.\Delta
= k$ and $k = - E_1.\Delta \geq 2$. The nefness of $L$ implies that $E.\Delta > 0$, whence from $2 = E.E_1 =
E.A + k E.\Delta \geq 2$ we deduce that $k = 2, E.\Delta = 1$ and $E.A = 0$. Hence $A \eqv qE$ for some $q
\geq 1$ by Lemma \ref{lemma10}. Now $2 = A.\Delta = q$ and therefore $E_1 \eqv 2E + 2\Delta$. We now set 
$E'_2 = E + \Delta$. Let us prove that $E + E'_2 = 2E + \Delta$ is nef. Let $\Gamma$ be a nodal curve such that
$(2E + \Delta).\Gamma < 0$. Since now $L \eqv 6E + 2\Delta$ the nefness of $L$ implies that $E.\Gamma > 0$. Now
$(2E + \Delta)^2 = 2$ and $(E + \Gamma)^2 \geq 0$ whence $(E + \Gamma).(2E + \Delta) \geq 1$. But this is a
contradiction since $(E + \Gamma).(2E + \Delta) = 1 + \Gamma.(2E + \Delta) \leq 0$.

\noindent Now that $E + E'_2$ is nef we just observe that by \cite[Prop.3.1.6 and Cor.3.1.4]{cd} either $E + E'_2$
or $E + E'_2 + K_S$ is base-component free, whence to conclude we choose accordingly $E_2 = E'_2$ or $E_2 = E'_2 +
K_S$.
\end{proof}
\renewcommand{\proofname}{Proof of Proposition {\rm \ref{g=9}}} 
 
\begin{proof}
We use the notation of Lemma \ref{lemmag=9}.

\noindent Suppose that there exists a line bundle $A$ that is a $g^1_4$ on $C$ and is such that $b_{2, A} =
0$. By the general remark \ref{genrem} we know that $C$ lies on a surface $X$ (obtained by desingularizing $Y_A$, if
necessary) and either $X = \Sigma_1, Bl_{V}Q$ and $C \sim -2K_X$ or $X = Bl_{V} C_8$ and $C$ is bielliptic. When
$X = \Sigma_1$ or $Bl_{V} C_8$ we get that $C$ has a complete base-point free $g^2_6$ and this is excluded by
\cite[Prop.3.5]{kl2}. The bielliptic case can also be excluded in another way, since, by \cite[Prop.4.17]{kl1}, $C$
has a unique $g^1_4$. Therefore $C \sim -2 K_X$ on $X = Bl_{V}Q$. By Proposition \ref{g1kondelp}(h) we have that $C$ has
a unique $g^1_4$, namely $f_{|C}$. Hence $(2E)_{|C} \sim f_{|C}$ and we deduce that $h^0((4E)_{|C}) = h^0(2f_{|C}) = 4$.
Now the exact sequence
\[ 0 \hpil - E_1 \hpil 4E \hpil (4E)_{|C} \hpil  0 \]
shows that $H^1(E_1 + K_S) \neq 0$, since $h^0(4E) = 3$. Therefore there exists a divisor $E_2$ as in Lemma 
\ref{lemmag=9}.

\noindent Let us prove that $(2E + E_2)_{|C}$ is a complete base-point free $g^2_8$ on $C$.

\noindent To this end note that since $(2E + E_2)^2 = 4$ and $2E + E_2$ is base-component free with two 
base points by Lemma \ref{lemmag=9} and \cite[Prop.3.1.6, Prop.3.1.4 and Thm.4.4.1]{cd}, we have that $(2E +
E_2)_{|C}$ is base-point free. Now the exact sequence
\[ 0 \hpil 2E + E_2 - C \hpil 2E + E_2 \hpil (2E + E_2)_{|C} \hpil  0 \]
shows that also $(2E + E_2)_{|C}$ is a complete $g^2_8$ since $2E + E_2 - C \eqv - 2 E - E_2$ whence 
$h^1(2E + E_2 - C) = 0$ because $2 E + E_2$ is nef.

\noindent Let $Z = E_2 \cap C$. Then $Z \subset C$ is an effective divisor such that 
$f_{|C} + Z \sim (2E + E_2)_{|C}$ is a complete base-point free $g^2_8$ on $C$, contradicting Proposition
\ref{g1kondelp}(j).
\end{proof}
\renewcommand{\proofname}{Proof}

\noindent We can now complete the proof of Theorem \ref{nodelp}. 

\renewcommand{\proofname}{Proof of Theorem {\rm \ref{nodelp}}}  
\begin{proof}
By Lemma \ref{lemmag=7}, Propositions \ref{g=7}, \ref{g^2_6particolari-tris}, \ref{g=8} and \ref{g=9} we can assume $L^2
\geq 18$.

\noindent Let $C$ be a curve as in the theorem, let $g = \frac{L^2}{2} + 1 \geq 10$ be the genus of $C$ and
suppose that $b_2(C) = 0$. By the general remark \ref{genrem} either $C$ is bielliptic or $g = 10$ and $C$ is
isomorphic to a smooth plane sextic. Now by \cite[Prop.4.17]{kl1} we have that $C$ has a unique $g^1_4$, therefore
it cannot be bielliptic. On the other hand the case of $C$ isomorphic to a smooth plane sextic is excluded
in \cite[Prop.3.1]{kl2}. Therefore we have a contradiction in all cases and the theorem is proved.
\end{proof}
\renewcommand{\proofname}{Proof}

\section{Proof of the main theorem}

\noindent We proceed with our main result. 

\begin{proof}
Let $C$ be a curve as in the theorem and let $g = \frac{L^2}{2} + 1 \geq 3$ be its genus.

\noindent Under the hypotheses (i) and (ii) the theorem follows immediately from Proposition \ref{gaussian}, while
if hypothesis (v) holds the theorem follows immediately from Corollary \ref{bel2}.

\noindent Now suppose we are under hypothesis (iii). By Theorem \ref{corkl1} and \cite[Prop.4.15]{kl1} we have that $C$ is
neither trigonal nor isomorphic to a smooth plane quintic, that is $\Cliff(C) \geq 2$. Then the theorem follows by
Proposition \ref{Cliff}.

\noindent Finally suppose that hypothesis (iv) holds. Since $L^2 \geq 12$, by \cite[Thm.1.4]{glm} (or by
Theorem \ref{corkl1}) we get that $\Cliff(C) \geq 2$. If $\Cliff(C) \geq 3$ then (iv) follows by Proposition
\ref{Cliff}(ii). If $\Cliff(C) = 2$ then, as is well-known, $C$ is either tetragonal or isomorphic to a smooth plane
sextic. But the latter case was excluded in \cite[Prop.3.1]{kl2}. Therefore $C$ is tetragonal and $\phi(L) = 2$ by Theorem
\ref{corkl1}. By Theorem \ref{nodelp} we have that $b_2(C) \geq 1$. Since $h^0(2K_C - M) = 1$ it follows that
$h^0(2K_C - M - b_2 A) = 0$ for every line bundle $A$ that is a $g^1_4$ on $C$. Therefore the theorem is a
consequence of Proposition \ref{tetragonal}(i).
\end{proof}


\begin{thebibliography}{[ACGH]}

\bibitem[ACGH]{acgh} E.~Arbarello, M.~Cornalba, P.~A.~Griffiths, J.~Harris. 
\textit{Geometry of Algebraic Curves, Volume I.}
Grundlehren der Mathematischen Wissenschaften \textbf{267}.
Springer-Verlag, New York, 1985.

\bibitem[Bd]{bd} L.~Badescu.
\textit{Polarized varieties with no deformations of negative weights}.
Geometry of complex projective varieties (Cetraro, 1990), Sem. Conf. \textbf{9}. Mediterranean,
Rende, 1993, 9-33.

\bibitem[BEL]{bel} A.~Bertram, L.~Ein, R.~Lazarsfeld. 
\textit{Surjectivity of Gaussian maps for line bundles of large degree on curves}. 
Algebraic geometry (Chicago, IL, 1989), 15--25, Lecture Notes in Math. \textbf{1479}. Springer,
Berlin, 1991.

\bibitem[BF]{bf} E.~Ballico, M.~Franciosi.
\textit{On Property $N_p$ for algebraic curves}. 
Kodai Math. J. \textbf{23}, (2000) 432--441.

\bibitem[BPV]{bpv} W.~Barth, C.~Peters, A.~van de Ven. 
\textit{Compact complex surfaces}.
Ergebnisse der Mathematik und ihrer Grenzgebiete \textbf{4}. 
Springer-Verlag, Berlin-New York, 1984. 

\bibitem[Bo]{bo} F.~Bogomolov. 
\textit{Holomorphic tensors and vector bundles on projective varieties}. 
Izv. Akad. Nauk SSSR Ser. Mat. \textbf{42}, (1978) 1227--1287, 1439.

\bibitem[C]{ca} F.~Catanese. 
\textit{Footnotes to a theorem of I. Reider}. 
Algebraic geometry (L'Aquila, 1988),  67--74. Lecture Notes in Math. \textbf{1417},
Springer, Berlin, 1990.

\bibitem[CD]{cd} F.~R.~Cossec, I.~V.~Dolgachev. 
\textit{Enriques Surfaces I}.
Progress in Mathematics \textbf{76}. Birkh\"auser Boston, MA, 1989.

\bibitem[CHM]{chm} C.~Ciliberto, J.~Harris, R.~Miranda. 
\textit{On the surjectivity of the Wahl map}.
Duke Math. J. \textbf{57}, (1988) 829--858.

\bibitem[CLM1]{clm1} C.~Ciliberto, A.~F.~Lopez, R.~Miranda. 
\textit{Projective degenerations of $K3$ surfaces, Gaussian maps, and Fano threefolds}.
Invent. Math. \textbf{114}, (1993) 641--667.

\bibitem[CLM2]{clm2}  C.~Ciliberto, A.~F.~Lopez, R.~Miranda.
\textit{Classification of varieties with canonical curve section via Gaussian maps on canonical
curves}.
Amer. J. Math. \textbf{120}, (1998) 1--21.

\bibitem[Co]{co} F.~R.~Cossec. 
\textit{On the Picard group of Enriques surfaces}.
Math. Ann. \textbf{271}, (1985) 577--600.

\bibitem[E]{ei} L.~Ein. 
\textit{The irreducibility of the Hilbert scheme of smooth space curves}. 
Proceedings of Symposia in Pure Math. \textbf{46}, (1987) 83--87.

\bibitem[F]{fr} R.~Friedman. 
\textit{Algebraic surfaces and holomorphic vector bundles}.
Universitext. Springer-Verlag, New York, 1998.

\bibitem[GH]{gh} P.~Griffiths, J.~Harris. 
\textit{Residues and zero-cycles on algebraic varieties}. 
Ann. of Math. \textbf{108}, (1978) 461--505.

\bibitem[GLM]{glm} L.~Giraldo, A.~F.~Lopez, R.~Mu\~{n}oz.
\textit{On the projective normality of Enriques surfaces}.
Math. Ann. \textbf{324}, (2002) 135--158.

\bibitem[Gr]{gr} M.~Green. 
\textit{Koszul cohomology and the geometry of projective varieties}. 
J. Differ. Geom. \textbf{19}, (1984) 125--171.

\bibitem[KL1]{klvan} A.~L.~Knutsen, A.~F.~Lopez. 
\textit{A sharp vanishing theorem for line bundles on K3 or Enriques surfaces}. 
Preprint 2005.

\bibitem[KL2]{kl1} A.~L.~Knutsen, A.~F.~Lopez. 
\textit{Brill-Noether theory of curves on Enriques surfaces I}. 
Preprint 2005.

\bibitem[KL3]{kl2} A.~L.~Knutsen, A.~F.~Lopez. 
\textit{Brill-Noether theory of curves on Enriques surfaces II}. 
In preparation.

\bibitem[KLM]{klm} A.~L.~Knutsen, A.~F.~Lopez, R.~Mu\~{n}oz. 
\textit{On the extendability of projective surfaces and a genus bound for Enriques-Fano
threefolds}.
Preprint 2006.

\bibitem[Kn1]{kn1} A.~L.~Knutsen. 
\textit{Exceptional curves on Del Pezzo surfaces}.
Math. Nachr. \textbf{256}, (2003) 58--81.

\bibitem[Kn2]{kn2} A.~L.~Knutsen. 
\textit{On $k$th order embeddings of $K3$ surfaces and Enriques surfaces}.
Manuscripta Math. \textbf{104}, (2001) 211--237. 

\bibitem [I1]{isk1} V.~A.~Iskovskih. 
\textit{Fano threefolds. I}.
Izv. Akad. Nauk SSSR Ser. Mat. \textbf{41}, (1977) 516--562, 717.

\bibitem [I2]{isk2} V.~A.~Iskovskih. 
\textit{Fano threefolds. II}.
Izv. Akad. Nauk SSSR Ser. Mat. \textbf{42}, (1978) 506--549.

\bibitem[L1]{la1} R.~Lazarsfeld.
\textit{A sampling of vector bundle techniques in the study of linear series}. 
Lectures on Riemann surfaces (Trieste, 1987). World Sci. Publishing, Teaneck, NJ, 1989,
500-559.

\bibitem[L2]{la2} R.~Lazarsfeld.
\textit{Lectures on linear series. With the assistance of Guillermo Fern\'andez del Busto}. 
IAS/Park City Math. Ser., 3, Complex algebraic geometry (Park City, UT, 1993),  161--219. 
Amer. Math. Soc., Providence, RI, 1997.

\bibitem[Lv]{lv} S.~L'vovsky. 
\textit{Extensions of projective varieties and deformations. I}.
Michigan Math. J. \textbf{39}, (1992) 41-51.

\bibitem[M]{mu} S.~Mukai. 
\textit{Biregular classification of Fano $3$-folds and Fano manifolds of coindex $3$}.  
Proc. Nat. Acad. Sci. U.S.A. \textbf{86}, (1989) 3000--3002.

\bibitem[N]{na} M.~Nagata. 
\textit{On rational surfaces. I. Irreducible curves of arithmetic genus $0$ or $1$}. 
Mem. Coll. Sci. Univ. Kyoto Ser. A Math. \textbf{32}, (1960) 351--370.

\bibitem[P]{pa} G.~Pareschi. 
\textit{Exceptional linear systems on curves on Del Pezzo surfaces}.
Math. Ann. \textbf{291}, (1991) 17--38.

\bibitem[R]{re} I.~Reider. 
\textit{Vector bundles of rank 2 and linear systems on algebraic surfaces}. 
Ann. of Math. \textbf{127}, (1988) 309--316.

\bibitem[S1]{sch1} F.~O.~Schreyer. 
\textit{Syzygies of canonical curves and special linear series}. 
Math. Ann. \textbf{275}, (1986) 105--137.

\bibitem[S2]{sch2} F.~O.~Schreyer. 
\textit{A standard basis approach to syzygies of canonical curves}. 
J. Reine Angew. Math. \textbf{421}, (1991) 83--123.

\bibitem[Te]{te} S.~Tendian. 
\textit{Gaussian maps on trigonal and tetragonal curves}. 
Unpublished.

\bibitem[T]{ty} A.~N.~Tyurin. 
\textit{Cycles, curves and vector bundles on an algebraic surface}. 
Duke Math. J. \textbf{54}, (1987) 1--26.

\bibitem[V]{vo} C.~Voisin. 
\textit{Courbes t\'etragonales et cohomologie de Koszul}. 
J. Reine Angew. Math. \textbf{387}, (1988) 111--121.

\bibitem[Wa]{wa} J.~Wahl. 
\textit{Introduction to Gaussian maps on an algebraic curve}.
Complex Projective Geometry, Trieste-Bergen 1989, London Math. Soc. Lecture Notes Ser. 
\textbf{179}. Cambridge Univ.\ Press, Cambridge 1992, 304-323. 

\bibitem[Za]{za} F.~L.~Zak. 
\textit{Some properties of dual varieties and their applications in projective geometry}.
Algebraic geometry (Chicago, IL, 1989), 273--280. Lecture Notes in Math. \textbf{479}. Springer,
Berlin, 1991.

\end{thebibliography}
\end{document}